\documentclass[journal]{IEEEtran}

\makeindex             

\usepackage{amsthm}
\usepackage{epsfig} 
\usepackage{epstopdf}
\usepackage{amsmath}
\usepackage{comment}
\usepackage{url}
\usepackage{algorithm,float}
\usepackage{algpseudocode}

\usepackage{amsfonts}
\usepackage{bm}
\usepackage{color}
\usepackage{cite}
\usepackage{stfloats}
\usepackage{chemarrow}
\usepackage{extarrows}

\usepackage{url}
\usepackage[switch]{lineno}
\usepackage{multirow}

\usepackage{setspace}

\DeclareGraphicsExtensions{.pdf,.png,.jpg,.eps,.ps}

\def\baa{\begin{align}}
	\def\eaa{\end{align}}

\newcommand{\bsq}{\begin{subequations}}
	\newcommand{\esq}{\end{subequations}}

\newcommand{\beq}{\begin{equation}}
	\newcommand{\eeq}{\end{equation}}
\newcommand{\bq}{\begin{eqnarray}}
	\newcommand{\eq}{\end{eqnarray}}
\newcommand{\bqn}{\begin{eqnarray*}}
	\newcommand{\eqn}{\end{eqnarray*}}
\newcommand{\bee}{\begin{enumerate}}
	\newcommand{\eee}{\end{enumerate}}
\newcommand{\bi}{\begin{itemize}}
	\newcommand{\ei}{\end{itemize}}

\newcommand{\diag}{\mathrm{diag}}

\newboolean{showcomments}
\setboolean{showcomments}{true}
\newcommand{\wang}[1]{\ifthenelse{\boolean{showcomments}}
	{ \textcolor[rgb]{1,0,1}{(ZW:  #1)}}{}}
\newcommand{\fliu}[1]{\ifthenelse{\boolean{showcomments}}
	{ \textcolor{red}{(FL:  #1)}}{}}
\newcommand{\zhao}[1]{\ifthenelse{\boolean{showcomments}}
	{ \textcolor{green}{(JP:  #1)}}{}}
\newcommand{\slow}[1]{\ifthenelse{\boolean{showcomments}}
	{ \textcolor{black}{(SL:  #1)}}{}}
\newcommand{\jpang}[1]{\ifthenelse{\boolean{showcomments}}
	{\textcolor[rgb]{0,0,0}{#1}}{}}
\newcommand{\jpangk}[1]{\ifthenelse{\boolean{showcomments}}
	{\textcolor[rgb]{0,0,0}{#1}}{}}
\newcommand{\jpangb}[1]{\ifthenelse{\boolean{showcomments}}
	{\textcolor[rgb]{0,0,0}{#1}}{}}
\theoremstyle{definition}

\theoremstyle{definition}

\newtheorem{remark}{Remark}

\newtheorem{example}{Example}

\ifCLASSOPTIONcaptionsoff
\usepackage[nomarkers]{endfloat}
\let\MYoriglatexcaption\caption
\renewcommand{\caption}[2][\relax]{\MYoriglatexcaption[#2]{#2}}
\fi
\allowdisplaybreaks
\begin{document}
	\graphicspath{{Paper_Fig/}}
	\setstretch{0.985}
	
	\title{Online Optimization in Power Systems with High Penetration of Renewable Generation: \\Advances and Prospects}

	 	\author{Zhaojian~Wang, 
		 		Wei Wei,
		 		John Zhen Fu Pang,
		 		Feng Liu,
		 		Bo Yang,
		 		Xinping Guan \IEEEmembership{Fellow,~IEEE}, \\and Shengwei Mei \IEEEmembership{Fellow,~IEEE}
		 }


	\maketitle
	
	\begin{abstract}   
		\jpang{Traditionally, offline optimization of power systems is acceptable due to the largely predictable loads and reliable generation. The increasing penetration of fluctuating renewable generation and Internet-of-Things devices allowing for fine-grained controllability of loads have led to the diminishing applicability of offline optimization in the power systems domain, and have redirected attention to online optimization methods. However, online optimization is a broad topic that can be applied in and motivated by different settings, operated on different time scales, and built on different theoretical foundations. This paper reviews the various types of online optimization techniques used in the power systems domain and aims to make clear the distinction between the most common techniques used. In particular, we introduce and compare four distinct techniques used covering the breadth of online optimization techniques used in the power systems domain, i.e., optimization-guided dynamic control, feedback optimization for single-period problems, Lyapunov-based optimization, and online convex optimization techniques for multi-period problems. Lastly, we recommend some potential future directions for online optimization in the power systems domain. }
	\end{abstract}
		
		\section*{Nomenclature}
		
		\begin{IEEEdescription}[\IEEEusemathlabelsep\IEEEsetlabelwidth{ADMM}]
			\item[AC] Alternating current
			\item[ADMM] Alternating direction method of multipliers
			\item[AGC] Automatic generation control
			\item[BFM] Branch flow model 
			\item[BIM] Bus injection model 
			\item[DC] Direct current
			\item[DG] Distributed generator
			\item[DR] Demand response
			\item[ESS] Energy storage system
			\item[GNG] Generalized Nash game
			\item[GNE] Generalized Nash equilibrium
			\item[KKT] Karush-Kuhn-Tucker
			\item[MPC] Model predictive control 
			\item[NE] Nash equilibrium
			\item[NEGDC] NE-guided dynamic control
			\item[OCO] Online convex optimization 
			\item[OGDC] Optimization-guided dynamic control 
			\item[OPF] Optimal power flow 
			\item[P2P] Peer-to-Peer
			\item[PV] Photovoltaic
			\item[RHC] Receding-horizon control 
			\item[RL] Reinforcement learning
			\item[SG] Synchronous generator 
			\item[SoC] State-of-charge 
			\item[WTG]  Wind power generator
		\end{IEEEdescription}

		\section{Introduction}
\begin{figure*}[!t]\label{Review_Structure}
			\centering
			\includegraphics[width=0.75\textwidth]{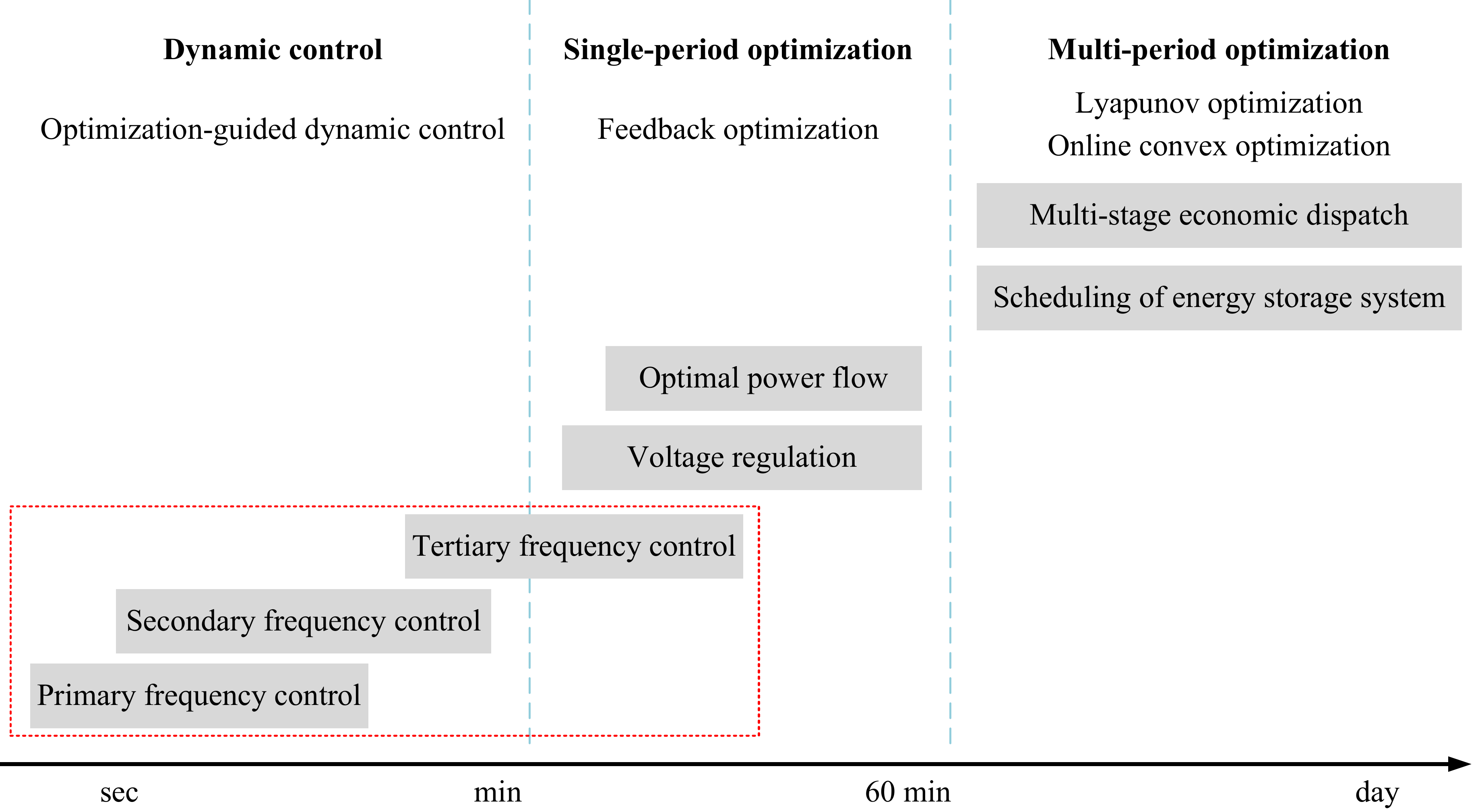}
			\caption{\color{black}Summary of online optimization in power systems. We find problems in the field of power systems can be grouped into three main classes of problems and highlight some archetypes here. For example, in the sub-minute time scale, solutions need to be obtained quickly, and therefore an optimization-guided dynamic controller is often used. On the other hand, problems that require the solving of high-dimensional, non-linear and complex power flow equations are often solved using feedback optimization. Lastly, problems that are solved at a higher time scale can usually be afforded a longer computation time and is suitable to be solved via Lyapunov optimization techniques or online convex optimization techniques. 
			}
		\end{figure*}

		\jpang{Global climate change has led to new global commitments to take active steps, i.e., COP21 Paris Agreement, on limiting global warming by at most two degrees celsius, compared to pre-industrial levels. As one of the primary drivers of climate change, there is immense pressure on countries to re-design their electric power systems to reduce their carbon footprint due to electricity generation by undertaking an increased penetration of renewable generation, e.g., wind and solar \cite{IRENA}. While these renewable generation sources may be operationally free of carbon emissions, either energy storage systems (ESSs) or alternative generating units that are highly controllable are needed to counteract the intrinsically volatile and uncertain generation. Exacerbating these conditions is the fact that these renewable generation sources are often integrated via power electronic devices with low inertia and rapid response speed. Recently technological advances in Internet-of-Things (IoT) devices have also proliferated a large number of high-power controllable loads, such as demand response and electric vehicles. Together, the effects brought about by the increase in renewable generation penetration and controllable loads present an unprecedented and severe operational challenge.}
		{\color{black}To this end, many \jpangb{(offline)} optimization methods have been \jpangb{proposed and} investigated, 
			\jpangb{which are largely} categorized into two types \cite{Zhongjie2021Real}: \jpangb{(i)} prediction-based and \jpangb{(ii)} historic data-based. For the first type, 
			forecasts of \jpangb{the renewable} generation are utilized in the problem formulation. Unfortunately, 
			\jpangb{as} the prediction \jpangb{of renewable generation} is not 
			accurate in practice, 
			\jpangb{these types of optimization methods} may lead to a sub-optimal solution. This 
			\jpangb{can result in poor} economic performance, \jpangb{and worst,} 
			can also cause stability issues in some cases. 
			For the second type, it mainly involves two kinds of approaches 
			\jpangb{distinguished by} the \jpangb{way uncertainty is modeled}
			from the historical data, namely, robust optimization and stochastic programming. In robust optimization, 
			uncertainty is modeled as a bounded set, which is determined in advance from the historical data. Then, the optimal solution \jpangb{which satisfies} 
			all of the possible cases \jpangb{within the bounded set} are selected. Consequently, the result is conservative and the suboptimality depends on the accuracy of the uncertainty set. In contrast, stochastic programming requires the distribution of uncertain variables, which is also obtained from historical data. The \jpangb{objective function} constraints including \jpangb{the corresponding} uncertainties are \jpangb{enforced to be} satisfied with a certain probability. However, the realization of uncertain variables may differentiate from the historical data, which leads to suboptimality or even infeasibility. \jpangb{As such, at times, it may be possible that} the
			results \jpangb{obtained from the offline optimization methods above} may not satisfy power flow equations or operational constraints in the 
			\jpangb{actual} situation. To summarize, 
			\jpangb{both} prediction-based and historic data-based methods \jpangb{face limitations in how they handle uncertainties, which may lead to suboptimality, infeasibility, or instability.}
			
			Recently, online optimization methods have attracted surging attention to overcome 
			\jpangb{the aforementioned} challenges, 
			\jpangb{as it} track changing conditions and thus are robust to uncertainties and variations \cite{molzahn2017survey}. In a power system with aggravating volatility and unpredictable uncertainties \jpangb{due to increasing penetration of renewable generation sources}, online optimization is required in many situations. To \jpangb{highlight the pervasiveness of}
			online optimization explicitly, we distinguish related works \jpangb{at} three \jpangb{different} time scale, 
			including 
			\jpangb{namely, }dynamic control, single-period problems, and multi-period problems. \jpang{Beyond the difference in time scale, these three types of problems are applied in different areas, and therefore have different constraints and objectives. These problems were initially solved using offline optimization techniques (in contrast to online optimization techniques), but the increase in renewable penetration has led to increased adoption of online optimization techniques to solve them.  }
			\jpang{The purpose of this review paper is to (i) highlight the different motivations and applications, (ii) bring to attention the different time scales, and (iii) present the theoretical foundations of these works. Here, we summarize them in Fig.\ref{Review_Structure}.}} 
		
		\textbf{1) Dynamic control}. 
		{\color{black}Frequency and voltage are two fundamental indices in power systems, which are \jpangb{typically} regulated by a hierarchical control structure \cite{Zhanqiang2022Predictive,Zhanqiang2022Event,Dorfler:Breaking}. 
		} In an alternating current (AC) power system, \jpang{frequency reveals the active power (im)balance} across the overall system, which varies from the nominal value as long as there is a power mismatch. A power system must maintain its frequency within a small neighborhood of the nominal value, typically \jpangk{50 or 60 Hz}. Otherwise, a power outage or even cascading failures in the power grid may happen \cite{maciver2021analysis}. \jpang{In the following, we take dynamic frequency control as an example. Traditionally, a hierarchical frequency control structure is adopted, consisting of three layers with different time scales, i.e., the primary control on several seconds time scale to stabilize the frequency, the secondary control on tens of seconds to a minute time scale to bring the frequency back to its nominal values, and the tertiary control from several minutes to tens of minutes \cite{jokic2007price} to solve for a more economical dispatch.
			Primary frequency regulation is traditionally performed using a simple proportional feedback controller, designed to limit the frequency deviation within an acceptable range rapidly. To eliminate frequency deviation, secondary frequency control is typically designed via automatic generation control (AGC). The first two layers intend to stabilize frequency and eliminate any deviation on a fast time scale, albeit potentially in a non-economical manner. Tertiary frequency control, often termed economic dispatch, focuses on an optimal generation allocation that addresses the power mismatch economically.
			In the past, this hierarchical structure works well due to the significant distinction in time scales for the different control schemes and the traditionally low uncertainty in the power system which meant that large frequency deviations are infrequent. 
			On the contrary, the increased volatility due to higher penetration of renewable energy sources in the grid requires optimality to be achieved at a much faster rate, thereby obfuscating the previously distinct hierarchy and time scale, which ultimately lead to the combination of different layers of control.} For example, tertiary control needs to be realized on the time scale of primary or secondary frequency control, which motivates the combination of the three to \jpangk{stabilize frequency, eliminate deviations,} and realize optimal operation simultaneously \cite{jokic2009constrained,Changhong:Design}. \jpangk{A majority of the work in this area can be described as an optimization-guided dynamic control (OGDC)}, also described as ``Breaking the hierarchy" \cite{Dorfler:Breaking} or ``Reverse engineering" \cite{Li:Connecting} in the field of power systems.

		\textbf{2) Single-period optimization}. \jpang{The bottleneck around many single-period problems like optimal power flow (OPF) and voltage regulation lies in solving the power flow equation. There are provably no analytical solutions for solving the power flow equation, and almost all algorithms in this field are offline and iterative in nature \cite{zohrizadeh2020survey}.
			Under these circumstances, it is necessary to wait for the converged result since the intermediate values obtained in these iterative processes in general do not satisfy the power flow equation, and thus, are inapplicable. The increased volatility due to the increase in renewable energy necessitates a time-efficient solution since the optimal solution in a previous timestep in a volatile system can be neither optimal nor feasible in a complex system like that of the power grid.} 
		\jpang{However, the power flow equation aims to mimic a high-dimensional, non-linear, and complex system, in which finding a solution is in general computationally costly.}
		\jpang{A recent idea rooted in optimization theory to overcome this dilemma is to exploit the laws of physics to solve the power flow equation in real-time and get the results not by optimization or solving an equation, but by simple-to-implement measurements.} This idea leads to the concept of feedback optimization \cite{gan2016online}, where measurements are collected, \jpangk{and used as part of} the algorithm design. 
		At each iteration, the algorithm computes generation set points and sends them to the corresponding generators. The grid then ``take these as input", and computes the states (generation, voltage, etc.) by the law of physics \jpangk{and electricity} in real-time. \jpang{The new measurements obtained from the grid will be used as the input in the next iteration. 
			A distinction from offline algorithms is that the intermediate iterates not only satisfy the power flow equations but can also be obtained at a fast time scale, which enables them to track evolving working conditions. }

		\textbf{3) Multi-period optimization}.
		The intrinsic volatility and uncertainty of \jpangk{the renewable sources of energy} call for sufficient backup resources, which lead to the proliferation of ESSs. The ESS scheduling problem is a multi-period problem with constraints like charging/discharging rates, the State-of-Charge of the battery (SoC), which can be a function of external conditions. For example, the use of district cooling systems as an energy storage system means that the charging and discharging rates, similar to the SoC, which is the current temperature of the facility, can be weather-dependent. As such, the decision variables are coupled in successive periods. {\color{black} Other constraints, like ramping constraints in the multi-stage economic dispatch, are also temporally coupled. Additionally, the action in any time step depends on variables in the future, such as power generation and price. 
			To this end, many approaches have been proposed, \jpangb{such as model predictive control which requires the use of prediction}, robust optimization, stochastic programming, and online optimization \cite{Zhongjie2021Real}. \jpangb{As introduced in the first paragraph of this section, results obtained from the first three methods are not fully amenable to uncertainties, and thus, suboptimality, infeasibility, and even instability may result. }
			On the other hand, online optimization allows decisions to be made after uncertain variables are observed or even without the information of uncertainties. In this way, the side effect of volatility and uncertainty is eliminated. Typical online optimization methods used to solve the multi-period optimization problem include Lyapunov optimization and online convex optimization (OCO), which will be introduced in detail in Section \ref{sec:mpo}.}
		
		\jpang{While online optimization is promising in the future power system with increased renewable penetration, the implications of its use still need to be further studied. As introduced above, the field of online optimization is broad and may seem inaccessible due to the different interpretations and massive amounts of corresponding works. To the best of our knowledge, there still lacks an explicit comparative analysis among these methods. In this paper, as shown in Fig.\ref{Review_Structure}, we introduce three main classes of problems in power systems that are amenable to online optimization techniques, and four different types of methods are presented comparatively in motivations, time scales, theoretic foundations, and typical applications. 
			The main merits of this paper are the following:}
		\jpang{
			\begin{itemize}
				\item We provide motivations for the use of online optimization within the power system and group online optimization problems into three main classes of problems, largely defined by their time scales, namely, dynamic control, single-period optimization, and multi-period optimization. Our hope is that it will eliminate the long-existing confusion in this field and lower the barrier of entry for researchers to work on online optimization in power systems. 
				\item We perform a comprehensive review and a comparative analysis of four types of online optimization methods, i.e., Optimization-guided Dynamic Control, Feedback Optimization, Lyapunov Optimization, and Online Convex Optimization, covering their motivation, time scale, theoretical foundations, and typical applications.
				\item We discuss critical challenges and promising directions for online decision-making in power systems in depth.
			\end{itemize}
		}
		
		\jpang{The rest of this paper is organized as follows.  
			Section II overviews background material on power systems detailing the power flow equations and power system dynamics, which aims to provide background and appreciation of the complexity of power flow equations to readers unfamiliar with the area of power systems. Section III reviews the main idea and typical applications of OGDC. Section IV summarizes the literature and techniques on feedback optimization for single-period problems, while Section V introduces online optimization methods for multi-period problems, including Lyapunov Optimization and OCO techniques. Lastly, we discuss the critical challenges and future directions of online optimization for power systems in Section VI, and conclude in Section VII.}
		

		\section{Preliminaries}
		\jpang{In this section, we introduce some preliminaries in power system optimization, including equations governing power flow and power system dynamics. This section serves two distinct purposes: (i) for completeness, we provide background information on power systems for those unfamiliar with the area for a richer appreciation of the complexity of power systems, and (ii) the provision of a concise set of equations that we will refer to in the sections to follow. }
		
		\subsection{Notations}
		
\begin{figure}[!t]
			\centering
			\includegraphics[width=0.45\textwidth]{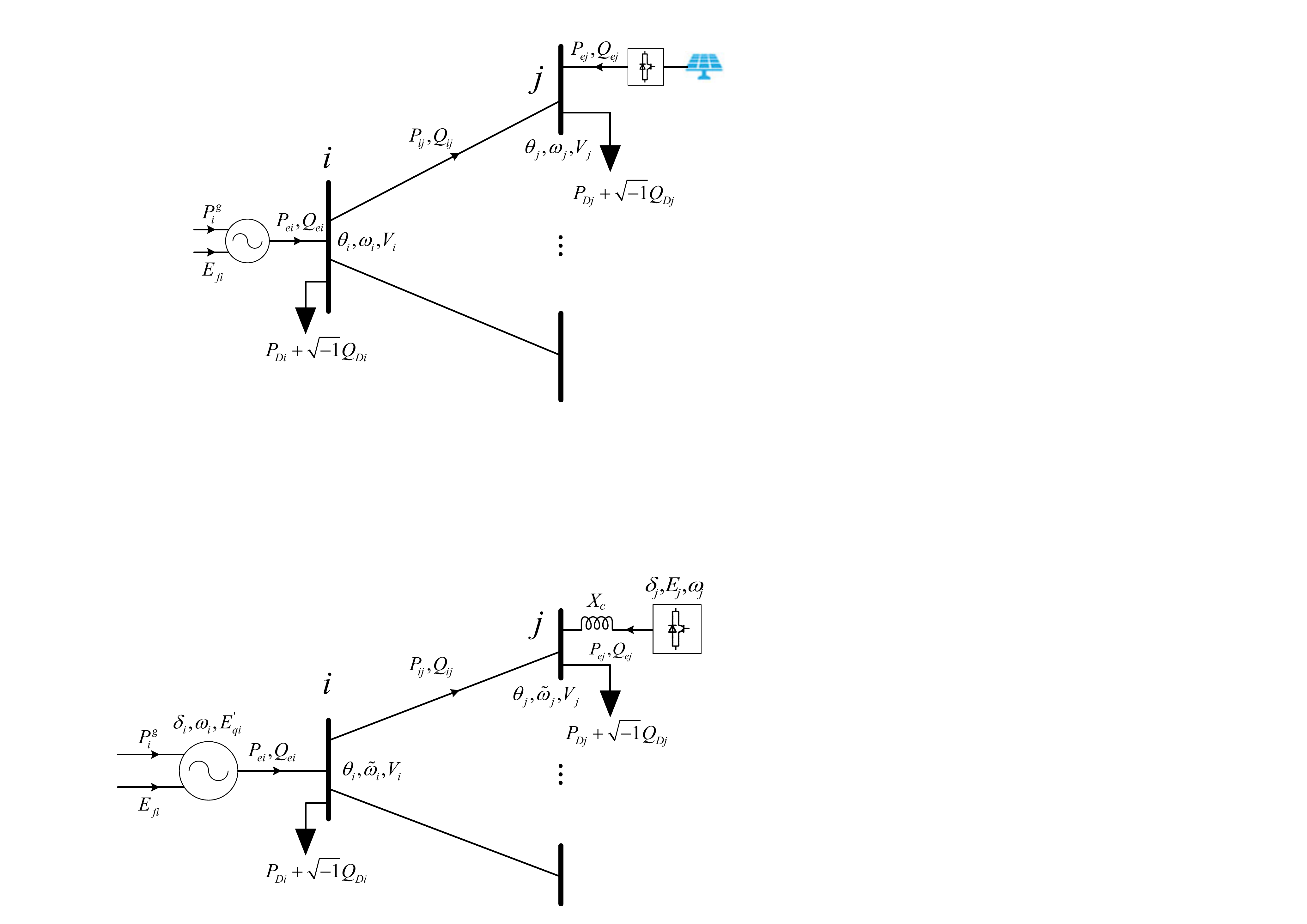}
			\caption{Summary of notations for a power system}
			\label{fig_notationsummary}
		\end{figure}

		\jpang{Consider an $n$-bus power system, succinctly illustrated in Fig.\ref{fig_notationsummary}. Each bus may consist of a combination of synchronous generators (SGs), inverter-integrated renewable generators, and loads. }
		Denote by $\mathcal{N}:=$ $\{1, \ldots, n\}$ the set of buses. Let $\mathcal E\subseteq \mathcal N\times \mathcal N$ be the set of lines, where line $(i,j)\in \mathcal E$ if buses $i$ and $j$ are connected directly. \jpang{Correspondingly, we use $\mathcal{N}_i$ to denote the set of buses that is connected to bus $i$.} Denote by $ r_{ij}+\mathbf{k}x_{ij} $ the impedance of line $ (i,j) $, where $\mathbf{k}:=\sqrt{-1}$. The network admittance matrix is denoted by $\mathbf{Y}:=\mathbf{G}+\mathbf{k B}$. The voltage at bus $i$ is $ {\bm V}_{i} := V_{i} e^{\mathbf{k} \theta_{i}} $, where $ V_{i} $ is the amplitude and $ \theta_{i} $ is the angle. Denote by $ {\bm I}_{ij} := I_{ij} e^{\mathbf{k} \theta_{ij}} $ the current over line $ (i,j) $, where $ I_{ij} $ is the amplitude and $ \theta_{i j} := \theta_{i}-\theta_{j} $ is the angle difference between $i, j$. We use $ P_{ij}, Q_{ij} $ to denote the active and reactive power from $i$ to $j$, whereas $\omega_i$ is the frequency, $P^g_i$ is the mechanical power input, $ E_{fi} $ is the excitation voltage, at bus $i$. The active and reactive power injected to bus $i$ by generators or inverters are denoted by $ P_{ei}, Q_{ei} $, respectively. The active and reactive loads attracted from bus $i$ are denoted by $ P_{Di}, Q_{Di} $, respectively. Then, the net power injections to the grid are $ P_i=P_{ei}-P_{Di},\ Q_i=Q_{ei}-Q_{Di} $. 
		Sometimes, we further distinguish active load into two types: controllable load and uncontrollable load, which are denoted by $ P_i^l, p_i $ with $ P_{Di}=P_i^l+p_i $, respectively. 
		
		\subsection{Power flow equations}\label{Power_flow_equations}
		\jpang{In this subsection, we introduce the power flow equations through both the bus injection model (BIM), branch flow model (BFM), and their linearized forms. }

		\subsubsection{Bus Injection Model }
		The power injection at the $ i $-th bus can be computed by
		%
		\begin{subequations}\label{Pre_BIM}
			\begin{align}\label{Pre_BIM1}
				P_{i} &=\sum\limits_{j\in N_i}  V_{i} V_{j} \left(G_{i j} \cos \theta_{i j}+B_{ij} \sin \theta_{i j}\right)\\
				\label{Pre_BIM2}
				Q_{i} &=\sum\limits_{j\in N_i}  V_{i} V_{j} \left(G_{i j} \sin \theta_{i j}-B_{ij} \cos \theta_{ij}\right)
			\end{align}
		\end{subequations}
		
		The corresponding active and reactive power flow $P_{ij},Q_{ij}$ from bus $i$ to bus $j$ are
		\begin{subequations}
			\label{line power}
			\begin{align}
				P_{ij}&= V_{i} V_{j} \left(G_{i j} \cos \theta_{i j}+B_{ij} \sin \theta_{i j}\right)
				\label{line_reactivepower}  \\ 
				Q_{ij}&=V_{i} V_{j} \left(G_{i j} \sin \theta_{i j}-B_{ij} \cos \theta_{ij}\right)
				\label{line_activepower}
			\end{align}
		\end{subequations}
		
		\jpang{Obviously, the power flow equations \eqref{Pre_BIM} and \eqref{line power} are nonlinear, which increase the complexity of many operation problems, as they introduce potential non-convexities which complicate optimization problems.
			To circumvent this, many works have contributed to linear approximations of the power flow equations above. These are made based on the following four key observations in high-voltage transmission networks:} 
		\begin{itemize}
			\item The reactive power over each line is very small compared to the active power counterpart, and thus could be omitted or neglected.
			\item The resistance is significantly less than the reactance, i.e., $r_{ij}\ll x_{ij}  $, which leads to the approximations $ G_{ij}\approx0$, and $\ B_{ij}=-\frac{1}{x_{ij}} $.
			\item For most  operating conditions, the difference in voltage angles of two buses $i,j,\forall (i,j)\in\mathcal{E}$ is very small, and therefore the linear function $x$ approximates the $\sin$ function well, i.e., $ \sin (\theta_{i}-\theta_{j})\approx\theta_{i}-\theta_{j} $.
			\item In the per-unit system, the voltage magnitude $V_i$ is very close to $1$ per unit.
		\end{itemize}
		\jpangk{Taking the above into consideration,} the active power flow over line $(i,j)$, i.e., $P_{ij}$, can be simplified into the following linearized form
		\begin{equation}\label{DC_line_reactivepower3}
			P_{ij}= B_{ij}\theta_{ij}
		\end{equation} 
		which is the so-called DC power flow. The DC power flow model is widely used in the optimal operation problems of high-voltage transmission networks. 
		
		\subsubsection{Branch flow model}
		As an alternative to BIM, balanced radial distribution networks can be represented using the BFM \cite{baran1989optimal}. Denote by $\ell_{i j}=I^2_{ij}$ the squared magnitude of the current flow from bus $i$ to bus $j$. \jpang{Under these scenarios, the branch flow equations can be written as }
		\begin{subequations}
			\label{DistFlow}
			\begin{align}
				\label{Active_Power}
				P_{ij}+P_j &=\sum\limits_{k \in  N_j }P_{jk}+r_{ij}\ell_{ij}\\
				\label{Reactive_Power}
				Q_{ij}+Q_j &=\sum\limits_{k \in  N_j }Q_{jk}+x_{ij}\ell_{ij} \\
				\label{V_Square}
				V_i^2-V_j^2 &=2(r_{ij}P_{ij}+x_{ij}Q_{ij})-\left(r_{ij}^2+x_{ij}^2\right)\ell_{ij}
				\\
				\label{S_Square}
				\ell_{ij}V_i^2&=P_{ij}^2+Q_{ij}^2
			\end{align}
		\end{subequations}
		
		\jpangk{The BFM is typically linearized by first} dividing two sides of \eqref{S_Square} by $V_i^2$, \jpangk{and substituting the resulting} $ \ell_{ij} $ into \eqref{Active_Power}, \eqref{Reactive_Power} and \eqref{V_Square} with the resulting equations, yielding $ r_{ij}\frac{P_{ij}^2+Q_{ij}^2}{V_i^2} $ and $x_{ij}\frac{P_{ij}^2+Q_{ij}^2}{V_i^2}$ as active and reactive power losses respectively. Their values are very small compared to line power flows, and hence are omitted hereafter. Defining a new variable $v_i:=\frac{V_i^2}{2}$, we obtain the following linearized BFM.
		\begin{subequations}
			\label{Linear_DistFlow}
			\begin{align}
				\label{Linear_Active_Power}
				P_{ij}+P_j&=\sum\limits_{k \in  N_j }P_{jk} \\
				\label{Linear_Reactive_Power}
				Q_{ij}+Q_j&=\sum\limits_{k \in  N_j }Q_{jk}  \\
				\label{Linear_V_Square}
				v_i-v_j &=r_{ij}P_{ij}+x_{ij}Q_{ij}
			\end{align}
		\end{subequations}
		As reported in \cite{baran1989optimal}, the approximation error  is usually on the order of $1\%$ . The linearized model \eqref{Linear_DistFlow} has been extensively used in the optimization and control of distribution networks due to its simplicity. 
		
		
		\subsection{Power system dynamics}
		Power system dynamics are determined by the types of equipment connected. With the proliferation of distributed generators (DGs) and demand response (DR), typical dynamic equipment include synchronous generators (SGs), inverters, and controllable loads. 
		\subsubsection{Synchronous generator}
		Here, we adopt the flux-decay model of the generator from \cite{Stegink:aunifying}, summarized in \eqref{Pre_eq:SGmodelC2}. 
		\begin{subequations}
			\label{Pre_eq:SGmodelC2}
			\begin{align}
				\dot \theta_i & =  \omega_i
				\label{Pre_eq:SGmodel.1aC2}
				\\
				M_i\dot \omega_i & =   P^g_i -D_i \omega_i- P_{Di} - \sum\nolimits_{j \in N_i} {P_{ij} } 
				\label{Pre_eq:SGmodel.1bC2}
				\\
				\begin{split}
					T_{d i}^{\prime} \dot{E}_{q i}^{\prime} &=E_{f i}-\left(1-\left(X_{d i}-X_{d i}^{\prime}\right) B_{i i}\right) E_{q i}^{\prime} \\
					&\quad-\left(X_{d i}-X_{d i}^{\prime}\right) \sum_{j \in \mathcal{N}_{i}} B_{i j} E_{q j}^{\prime} \cos \theta_{i j}
				\end{split}
				\label{Pre_eq:SGmodel.1cC2} 
				\\
				\label{Pre_eq:SGmodel.1d}
				T^g_i\dot P_i^g &= - P_i^g + u_i^g - {\omega_i(t)}/{R_i}
			\end{align}
		\end{subequations}
		Here, $M_i$ is the moment of inertia, $D_i$ is the damping constant,  $T^{'}_{di}$ is the $d$-axis transient time constant, and \eqref{Pre_eq:SGmodel.1aC2}, \eqref{Pre_eq:SGmodel.1bC2} are the so-called swing equations. Equation \eqref{Pre_eq:SGmodel.1d} is the simplified model of the governor and turbine. $ X_{di} $ is the $d$-axis synchronous reactance, and $ X_{di}^{\prime} $ is the $d$-axis transient reactance.
		
		\subsubsection{Inverter}
		The dynamics of power inverters are determined by the control strategies adopted. Droop control is the most common controller and is widely used in power inverters \cite{Schiffer:Conditions}. Based on droop control, the \jpangk{dynamics} at bus $j$ is
		\begin{subequations}
			\label{DG_dynamic2}
			\begin{align}
				\label{DG_dynamic2_delta}
				{\dot \theta _j}  &=  {\omega _j}
				\\
				\label{DG_dynamic2_omega}
				{\tau _j}{{\dot \omega }_j}  &= -k_{{\omega_j}}\omega_j-\sum\nolimits_{k \in N_j} {P_{jk} }-P_{{Dj}} +u_{j}^P
				\\
				\label{DG_dynamic2_E}
				{\tau _j}{\dot V_j} & =  -k_{{V_j}}V_j-\sum\nolimits_{k \in N_j} {Q_{jk} } -Q_{{Dj}} + u_{j}^Q
			\end{align}
		\end{subequations}
		where ${\tau _j}, k_{{\omega_j}}, k_{{V_j}}$ are positive constants, and $u_{j}^P, u_{j}^Q$ are control inputs. 
		
		{\color{black}Besides \jpangb{the traditional} droop control, there are many control strategies for inverter-based DGs, such as quadratic droop control, reactive current control, and exponentially scaled averaging reactive power control. As they are not commonly adopted as droop control, we do not introduce them in detail. An interested reader can refer to \cite{de2018bregman} for details. }
		\subsubsection{Load bus}
		For a pure load bus, if the controllable load is considered, its dynamics could be designed as needed, such as via an inertia link as in \cite{wang2019distributed1}. 
		\begin{align}
			T^l_j \dot P^{l}_j & =  - P^{l}_j(t) + u^l_j
		\end{align}
		where $P_j^l$ is the controllable load, $u_j^l$ is the control input, and $T^l_j$ is the time constant of the controllable load. Uncontrollable loads are usually simply modeled as a constant.

		\section{Optimization-guided Dynamic Control}\label{model}
		\jpang{In this section, we first introduce the main idea of Optimization-guided Dynamic Control (OGDC). Then, we survey applications of OGDC in the power systems area. }
		\subsection{Main idea}\label{Idea_realtime}
\begin{figure}[!t]
			\centering
			\includegraphics[width=0.42\textwidth]{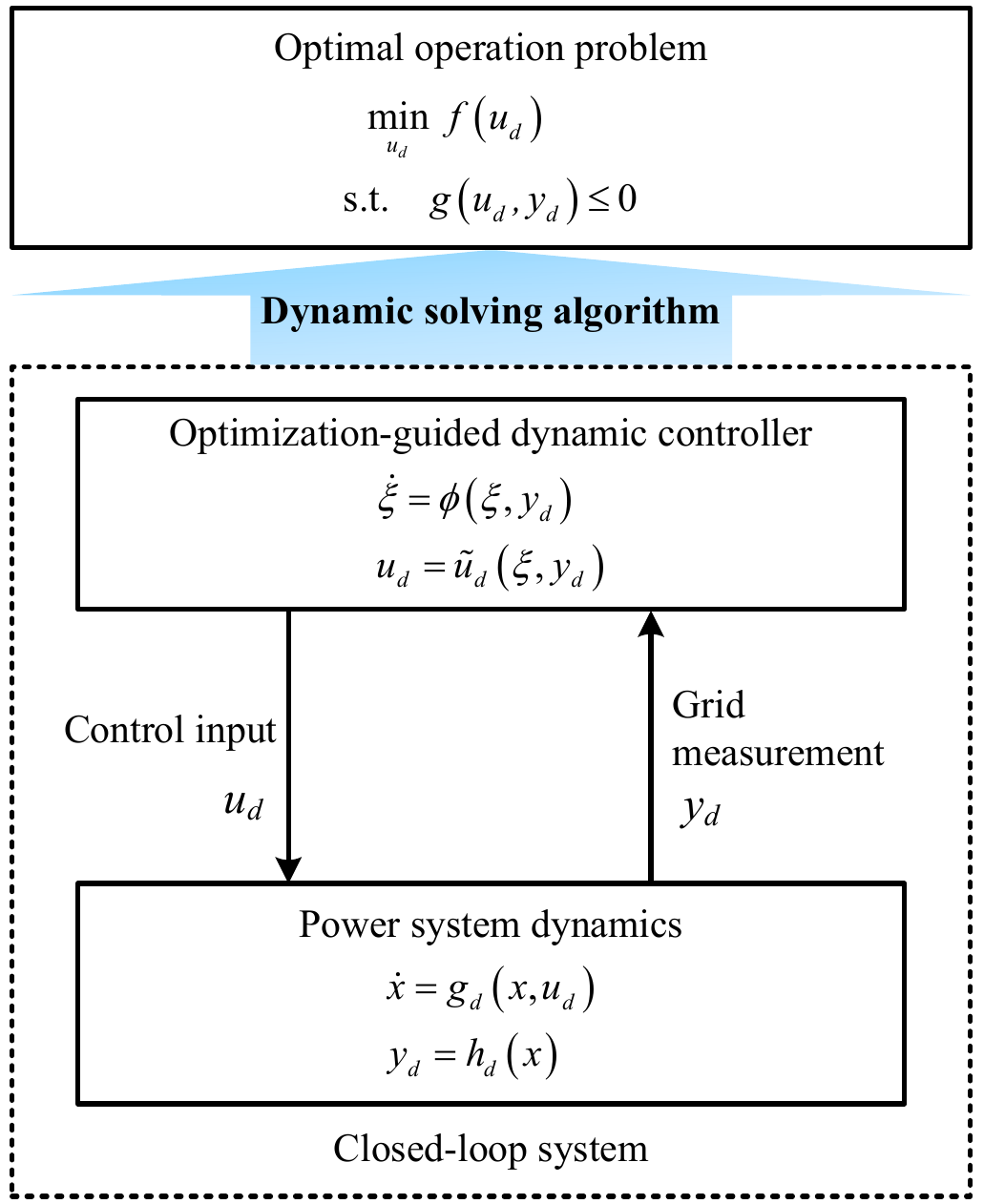}
			\caption{Diagram of optimization-guided dynamic control}
			\label{Feedback_control}
		\end{figure}

		\jpang{The main idea of optimization-guided dynamic control is to  design a (dynamic) feedback controller for physical systems, which will drive the system states based on the optimal solution of an optimization problem. The overall framework of OGDC is illustrated in Fig.\ref{Feedback_control}.}
		
		In this framework, the lower layer is the fast time-scale dynamics of a physical system. Its state variables ($ x $), control inputs ($ u_d $), and outputs ($ y_d $) are determined by the integrated devices. With reference to the dynamics introduced in Section II.C, the relevant control inputs are $ u_i^g, u_{j}^P, u_{j}^Q, u^l_j$. The output measurements $ y_d $  could correspondingly include frequency, voltage, and active and reactive power. 
		The upper layer is the optimal operation problem, such as tertiary frequency control.

		The general form of the optimization problem is 
		\begin{subequations}
			\label{Optimiztion_guided_control_optimization_problem}
			\begin{align}
				\mathop {\min }\limits_{{{ u}_d}}\quad & f\left( {{{ u}_d}} \right)\\
				\text{s.t.}\quad & g\left( {{{ u}_d},{{ y}_d}} \right) \le 0
			\end{align}
		\end{subequations}
		where the objective function is to minimize the control cost, and the constraint is associated with the input and output of the system. 
		
		Then, a dynamic feedback controller is designed:
		\begin{align}\label{eq:Intro_equivalent_controller}
			\left\{
			\begin{array}{rcl}
				\dot {\xi}  &=& \phi(\xi, y_d) \\
				u_d &=& \tilde{u}_d(\xi, y_d)
			\end{array}
			\right. 
		\end{align}
		By applying the dynamic feedback controller to the physical system, we can get the closed-loop system
		\begin{subequations}\label{eq_closedloop}
			\begin{align}
				\dot{x}&=g_{d}\left(x, \tilde{u}_d(\xi, y_d)\right)\\
				\dot {\xi}  &= \phi(\xi, y_d) \\
				y_{d}&=h_{d}(x)
			\end{align}
		\end{subequations}
		
		Our goal is to steer the closed-loop system \eqref{eq_closedloop} to the optimal solution of the optimization problem \eqref{Optimiztion_guided_control_optimization_problem} in the steady state. The equivalent interpretation is that \eqref{eq_closedloop} defines a dynamic algorithm which solves \eqref{Optimiztion_guided_control_optimization_problem} implicitly. 
		
		To the best of our knowledge, this idea dated back to \cite{jokic2009constrained} in the field of power systems, which presented a methodology to regulate the system to an optimal operation point, i.e., a solution to a given economic dispatch problem. Then, \cite{jokic:real} further exploits the shadow price interpretation of Lagrange multipliers to guarantee the economically optimal operation of power systems. This idea has recently been widely used in optimal frequency control. 
		We use the optimal secondary frequency control in \cite{wang2019distributed1} as an example to illustrate how to design an optimization-guided dynamic controller. 
		
		\begin{example}
			{\color{black}Increasing the penetration of volatile renewable generations causes power imbalance to fluctuate rapidly with a large magnitude, requiring tertiary frequency control to be realized on a fast time scale, coinciding with that of primary or secondary frequency regulation. }
			To design the controller, the DC power flow model introduced in \eqref{DC_line_reactivepower3} is adopted, where voltage dynamics is not considered. Then, the frequency dynamics with SG and controllable loads are simplified into the following form:
			\begin{subequations}\label{eq:model.1}
				\begin{align}
					\dot \theta_j & =  \omega_j
					\label{eq:model.1a}
					\\
					\begin{split}
						M_j \dot \omega_j & =   P^g_j- P^{l}_j - p_j -D_j \omega_j
						\\
						& \quad + \sum_{i: i\rightarrow j} B_{ij}\theta_{ij}
						-  \sum_{k: j\rightarrow k}  B_{jk}\theta_{jk}
					\end{split}
					\label{eq:model.1b}
					\\
					T^g_j \dot P^g_j & =  - P^g_j + u^g_j - {\omega_j(t)}/{R_j}
					\label{eq:model.1c}
					\\
					T^l_j \dot P^{l}_j & =  - P^{l}_j(t) + u^l_j
					\label{eq:model.1d}
				\end{align}
			\end{subequations}

			As discussed earlier, we first define the control goal as an optimization problem, as in \eqref{Optimiztion_guided_control_optimization_problem}. 
			\begin{subequations}\label{eq:opt.1} 
				\begin{align}
					\min & \quad \frac{1}{2} \sum\limits_{j\in \mathcal{N}} \alpha_j \left(P^g_j\right)^2 
					+ \frac{1}{2} \sum\limits_{j\in \mathcal{N}} \beta_j \left(P^l_j\right)^2
					+ \frac{1}{2} \sum\limits_{j\in \mathcal{N}} D_j \omega_j^2
					\nonumber
					\\
					\label{eq:opt.1d}
					\\
					\text{over} &\quad {\theta}, \omega, P^g, P^l
					\nonumber
					\\ 
					\text{s. t.}  
					& \label{eq_local_balance}\quad P^g_j =  P^l_j + p_j, \qquad j\in \mathcal{N}\\
					\label{eq:opt.1a} 
					\begin{split}
						&\quad P^g_j = P^l_j  + p_j +  D_j \omega_j   \\
						& \qquad\quad - \sum\limits_{i: i\rightarrow j} B_{ij} \theta_{ij}
						+ \sum\limits_{k: j\rightarrow k} B_{jk} \theta_{jk}, \ j\in \mathcal{N}
					\end{split}
					\\
					&\quad \underline{P}^g_j  \leq \ P^g_j(t) \ \leq \overline{P}^g_j, \quad j\in \mathcal{N}
					\label{eq:OpConstraints.1a}\\
					&\label{eq:OpConstraints.1b}\quad \underline{P}^l_j  \leq \ P^l_j(t) \ \leq \overline{P}^l_j, \quad j\in \mathcal{N}
					\\
					& \quad P^g_j \ = \ u^g_j, \quad j\in \mathcal{N}
					\label{eq:opt.1b}
					\\
					& \quad P^l_j \ = \ u^l_j, \quad j\in \mathcal{N}
					\label{eq:opt.1c}
				\end{align}     
			\end{subequations}
			where $\alpha_j>0$, $\beta_j>0$ are constant weights, $\underline{P}_{j}^{g}, \bar{P}_{j}^{g} $ are lower and upper bounds of $P_{j}^{g}$, and $\underline{P}_{j}^{l}, \bar{P}_{j}^{l} $ are lower and upper bounds of $P_{j}^{l}$.
			The first two terms in the objective function aim to minimize the regulation cost of generators and controllable load, while the last term is part of the frequency controller design. The constraint \eqref{eq_local_balance} is the local power balance, and \eqref{eq:opt.1a} is obtained from \eqref{eq:model.1b} by setting $\dot\omega_j=0$.
			Here, \eqref{eq:OpConstraints.1a} and \eqref{eq:OpConstraints.1b} are hard limits on the regulation capacities of generation and controllable load at each node, which should not be violated at any time even during the transient period. To ensure that these hard constraints are not violated, any violations will be projected onto the desired region. Lastly, constraints \eqref{eq:opt.1b}, \eqref{eq:opt.1c} reveal the relationship between control input and states.  
			
			For each node $j\in \mathcal{N}$, the control law is
			\begin{subequations}\label{Optimization_controller}
				\begin{align}
					\dot \xi_j & =  \gamma^{\xi}_j \left( P^g_j -  P^l_j - p_j\right)
					\label{eq:control.1c}
					\\
					u^g_j & =  \left[ 
					P^g_j - \gamma^g_j \left( \alpha_j P^g_j + \omega_j + \xi_j \right)
					\right]_{\underline P^g_j}^{\overline P^g_j}  + \frac{\omega_j}{R_j}
					\label{eq:control.1a}
					\\
					u^l_j & =  \left[ 
					P^l_j - \gamma^l_j \left( \beta_j P^l_j - \omega_j - \xi_j \right)
					\right]_{\underline P^l_j}^{\overline P^l_j}
					\label{eq:control.1b}
				\end{align}
			\end{subequations}
			where $\gamma^g_j, \gamma^l_j, \gamma^{\xi}_j$ are positive constants. 
			For any $x_i, a_i, b_i \in \mathbb R$ with $a_i\leq b_i$, we define the operator
			$$[x_i]_{a_i}^{b_i} := \min \{ b_i, \max \{ a_i, x_i \} \}$$ 
			
			If we further set the control gains as 
			$ \gamma^g_j \ = \ ( T^g_j )^{-1}, 
			\gamma^l_j \ = \ ( T^l_j )^{-1}
			$ and apply \eqref{Optimization_controller} to the dynamic system \eqref{eq:model.1}, the closed-loop system is (in vector form)
			\begin{subequations}\label{eq:model.2}
				\begin{align}
					\dot {\tilde\theta}&=  C^T\omega \label{eq:model.2a'}\\
					\dot \omega &= M^{-1}\left(P^g-P^l-p-D\omega(t)-CB\tilde\theta \right )  \label{eq:model.2b2'}	\\
					\dot{P}^g&= (T^g)^{-1}\left ( -P^g+\hat u^g \right)\\ 
					\dot{P}^l&= (T^l)^{-1}\left ( -P^l+\hat u^l \right)\\
					\dot{\xi}&= \Gamma^{\xi} \left (P^g-P^l-p \right)
				\end{align}
			\end{subequations}
			where 
			\begin{align*}
				\hat u^g & =  \left[ P^g - (T^g)^{-1} \left( 
				A^g P^g + \omega + \xi \right)\right]_{\underline P^g}^{\overline P^g}
				\\
				\hat u^l & =  \left[ P^l - (T^l)^{-1} \left( 
				A^l P^l - \omega - \xi \right)\right]_{\underline P^l}^{\overline P^l}
			\end{align*}
			and $M := \diag(M_j, j\in N)$, $ T^l:=\diag(T^l_j, j\in N )$,  $T^g:=\diag(T^g_j, j\in N ),\Gamma^{\xi} := \diag(\gamma^{\xi}_j, j\in N)$, $A^g := \diag(\alpha_j, j\in N)$, $A^l := \diag(\beta_j, j\in N)$, and $B:=\diag(B_{ij}, (i,j)\in E)$. 
			
			\jpang{It is shown in \cite{wang2019distributed1} that the closed-loop system \eqref{eq:model.2} serves as the partial primal-dual gradient dynamics which solves the original problem \eqref{eq:opt.1}. The asymptotic stability of \eqref{eq:model.2} can  also be proved theoretically. Thus, the controller \eqref{Optimization_controller} can indeed drive the physical system to an optimal operating state. This is why such controllers are known as ``optimization-guided dynamic control".}
			
		\end{example}

		\subsection{Application survey}
		\begin{table*}[!t]
			\small
			\renewcommand{\arraystretch}{1.3}
			\centering
			\caption{\textsc{Comparison of problem setups}}
			\label{Comparison_OGC}
			\begin{tabular}{c c c}
				\hline
				\hline 
				& & References \\
				\hline
				\multirow {2}{*}{Hard input limit}  & Yes & \cite{Changhong:Design,Mallada-2017-OLC-TAC,Kasis:Primary1,kasis2021primary,wang2019distributed1,kasis2019stability,kasis2020secondary,wang2020distributednonsmooth} \\
				& No & \cite{Li:Connecting,wang2019distributed,wang2019distributed_variation,trip2016internal,chen2016reverse,Stegink:aunifying,de2019feedback} \\
				\hline
				\multirow {2}{*}{Line congestion}  & Yes & \cite{Mallada-2017-OLC-TAC,wang2019distributed1,wang2020distributednonsmooth,Stegink:aunifying} \\
				& No & \cite{Changhong:Design,Kasis:Primary1,kasis2021primary,kasis2019stability,kasis2020secondary,Li:Connecting,wang2019distributed,wang2019distributed_variation,trip2016internal,chen2016reverse,de2019feedback}\\
				\hline
				\multirow {2}{*}{Load type}  & Constant & \cite{Mallada-2017-OLC-TAC,wang2019distributed1,wang2020distributednonsmooth,Stegink:aunifying,Changhong:Design,Kasis:Primary1,kasis2021primary,kasis2019stability,kasis2020secondary,Li:Connecting,wang2019distributed,chen2016reverse,de2019feedback} \\
				& Time-varying & \cite{wang2019distributed_variation,trip2016internal} \\
				\hline
				\multirow {2}{*}{Structure model}  & Reduced & \cite{wang2019distributed_variation,trip2016internal,wang2019distributed1,Stegink:aunifying,kasis2021primary,kasis2020secondary,Li:Connecting,de2019feedback} \\
				& Preserved & \cite{Mallada-2017-OLC-TAC,wang2020distributednonsmooth,Changhong:Design,Kasis:Primary1,kasis2019stability,wang2019distributed,chen2016reverse} \\
				\hline
				\multirow {2}{*}{System model}  & Linear & \cite{wang2019distributed_variation,Mallada-2017-OLC-TAC,wang2019distributed1,wang2020distributednonsmooth,Changhong:Design,kasis2021primary,kasis2019stability,Li:Connecting,chen2016reverse} \\
				& Nonlinear & \cite{trip2016internal,Stegink:aunifying,Kasis:Primary1,kasis2020secondary,wang2019distributed,de2019feedback} \\
				\hline
			\end{tabular}
		\end{table*}

		In this subsection, we will introduce more applications of optimization-guided control in the frequency regulation of AC systems. 
		We present here a survey on two recent methods used to design OGDC in power systems: consensus methods and primal-dual gradient methods. 
		
		\subsubsection{Optimal frequency control based on consensus method}

		In consensus-based control, the agents, such as generators, loads, MGs, or other forms of local systems, estimate a global variable using a consensus algorithm \cite{Olfati:Consensus}. Specifically, in power systems, if we take the marginal cost as the global variable, all agents share the same marginal cost in the steady state\footnote{Another global variable is the ratio between actual generation and the maximal capacity, which implies that all generators supply the load fairly up to their maximal capability. Because it is not optimization-guided, we do not introduce it in this paper. Readers could refer to  \cite{Xin:A, Guo2015Distributed, Simpson2015Secondary} for more details. }. This implies that the generation configuration is economically optimal. 
		
		A consensus-based controller is designed in \cite{Wu2018A}, which takes the following form. 
		\begin{subequations}
			\begin{align}\label{eq_consensus}
				\dot {\xi}_i  &= -c\left[\sum\nolimits_{j \in N_{i}} a_{i j}\left(\omega_{i}-\omega_{j}\right)+\left(\omega_{i}-\omega_{r e f}\right)\right] \nonumber\\ 
				&\qquad-k \sum\nolimits_{j \in N_{i}} a_{i j}\left( \eta_{i}\left(P_{i}\right)- \eta_{j}\left(P_{j}\right)\right)\\
				u^P_i&=\xi_i
			\end{align}
		\end{subequations}
		where $ c $ and $ k $ are constant, $ a_{ij} $ is the weight between bus $i$ and $j$, $ \eta_{i}\left(P_{i}\right) $ is the marginal cost of agent $i$, and $\omega_{ref}$ is the frequency reference. {\color{black}Under this setting, the system frequency will be restored to the rated value $\omega_{ref}$, and equal marginal cost is obtained, which realizes both frequency stability and economical allocation of active power among all agents. In \cite{Changhong2015ACC}, the control input $u_i^g$ is defined as the integral of marginal cost differences and frequency deviations, i.e., $u_i^g=-\int_0^t \omega_i(\tau) d\tau+\int_0^t\sum\nolimits_{j\in{N_i}}(\eta_i(\tau)-\eta_j(\tau))d\tau $ with $ \eta_i $ as the marginal cost. With that, the nominal system frequency is restored and the marginal cost reaches a consensus among all of the participating generators.} Similar ideas are also adopted in \cite{cady2015distributed,Dorfler:Breaking,dorfler2017gather,weitenberg2019robust} for optimal frequency control, and further improved in \cite{Simpson2015Secondary,Wu2018A,Watson2021Frequency} by considering voltage control simultaneously, and in \cite{wang2020asynchronous,lai2020stochastic,cintuglu2021real} considering asynchronous information.
		
		The consensus-based optimal frequency control can achieve an identical marginal cost among all agents, which, however, also limits its applications. For example, the marginal cost will not be identical if line power limits exist.  {\color{black}In this situation, the consensus-based method will not work, but it motivates primal-dual gradient methods, which we present next.}

		\subsubsection{Optimal frequency control based on primal-dual gradient method}
		
		The main idea is to use the primal-dual gradient dynamics, i.e., the well-known saddle dynamics \cite{Cherukuri:Asymptotic,colombino2019online,wang2021exponential,Xinlei2022A}, to solve the optimization problem \eqref{Optimiztion_guided_control_optimization_problem}. Then, one can combine the solution dynamics with the power system dynamics, and the closed-loop system converges to an equilibrium that corresponds to the optimal solution of the original optimization problem \cite{Changhong:Design, zhang2018distributed,Li:Connecting}. In this way, it realizes optimization-guided dynamic control.
		
		As introduced in Section III-A, this idea takes root in \cite{jokic2009constrained,jokic:real}, which uses the Karush-Kuhn-Tucker (KKT) condition to regulate a nonlinear dynamical system to the optimal solution of a given optimization problem. This idea is further generalized by using the primal-dual gradient dynamics. A rich literature has emerged investigating optimal frequency control based on primal-dual gradient method, which can be roughly divided into two categories: primary-tertiary control \cite{Changhong:Design,Mallada-2017-OLC-TAC,Kasis:Primary1,kasis2021primary} and secondary-tertiary control \cite{trip2016internal,chen2016reverse,Li:Connecting,wang2019distributed1,wang2019distributed,kasis2019stability,wang2019distributed_variation,kasis2020secondary,wang2020distributednonsmooth,Stegink:aunifying,de2019feedback}. We present a detailed classification of these works in Table \ref{Comparison_OGC} according to indices including input limits, line congestion, load type, structure model, and system model. 
		
		In the first category, the primary and tertiary frequency control is combined, which intends to stabilize the system rapidly with the best economic efficiency. In \cite{Changhong:Design}, an optimal load-side control problem is formulated and a dynamic controller based on the partial primal-dual gradient method is derived, which provides a paradigm shift for such research. It is extended in \cite{Mallada-2017-OLC-TAC} to consider the line congestion, in \cite{Kasis:Primary1} to relax the model requirements, and in \cite{kasis2021primary} to consider on-off loads. Since the response speed is the most important task of the primary frequency control, the primary-tertiary control is usually decentralized without the need for communication. The main trade-off of such an approach is that frequency deviation persists. 
		
		In the second category, the secondary and tertiary frequency control is combined, with the intention to both restore nominal frequency as well as realize economic dispatch. In \cite{Li:Connecting}, the notion of reverse engineering is proposed, which interprets primal-dual gradient dynamics as solving economic dispatch problems to the power system dynamics together with AGC. Consequently, the secondary frequency control can achieve economic dispatch simultaneously. This idea is extended in a similar fashion to settings with various practical considerations, such as hard operation constraints in \cite{wang2019distributed1}, partial control coverage in \cite{pang2017optimal,wang2019distributed}, more general physical dynamic models \cite{kasis2019stability}. Further, the time-varying disturbances could also be considered in this framework \cite{trip2016internal,wang2019distributed_variation}. Some supplementary controllers can be added to deal with the disturbances, such as the internal model control \cite{trip2016internal,wang2019distributed_variation} and high-gain observer \cite{CHOWDHURY2021109753}. {\color{black}Besides the economic dispatch, the power market and secondary frequency control also can be realized simultaneously, where the power market is solved by a dynamic algorithm  \cite{Stegink:aunifying,de2019feedback,cherukuri2021frequency}, where \cite{de2019feedback} considers a Cournot competition market model and \cite{cherukuri2021frequency} a Bertrand competition market model. }
		
		Primal-dual gradient dynamics are very appealing since they can partially be interpreted as system physics. It also has wider applicability in power systems, which is also utilized in the DC systems, including the unified OPF control in \cite{wang2017unified} and emergency control in hybrid AC-DC grids  \cite{Ye2021Optimal}. 
		
		\begin{remark}
			OGDC focuses on the optimality in the steady state on the time scale of dynamic control, which is complementary to the traditional optimal control, such as the linear quadratic regulator. The latter intends to optimize the integral of a quadratic form, such as 
			\begin{align}\label{transient_performance}
				\min\quad\frac{1}{2} \int_{0}^{\infty}\left( x^{T}(t)R_xx(t)+ u^{T}(t)R_uu(t)\right) \mathrm{d} t
			\end{align}
			with positive definite $ R_x,R_u $.  
			The objective \eqref{transient_performance} is to minimize the control cost over the whole transient process. 
			These two aspects are both very important when large amounts of volatile renewable generations are integrated. However, it is still an open question of how to achieve optimal performance on both sides simultaneously. 
		\end{remark}

		\section{Feedback Optimization for Single-Period Problems}\label{Properties}
		In this section, we first introduce the main idea of feedback optimization for single-period problems before surveying applications in OPF and voltage regulation. 
		
		\subsection{Main idea}
		Feedback optimization is an alternative way to solve single-period optimization problems in time-varying environments, especially for problems which contain constraints of power flow equations. 
		This is because this algorithm is able to solve power flow equations on a fast time scale, which is required for the power imbalance in the system caused by volatile loads and rapidly fluctuating renewable generation sources. It is widely known that power flow equations are nonlinear and high-dimensional, and therefore usually computationally costly. 
		In this situation, traditional offline approaches may not be applicable. \jpang{Fortunately, a physical power system by itself turns out to be a very efficient computer capable of calculating the exact solutions to power flow equations. Thus, one can get the results directly by measurement, which is more effective and accurate than computing or calculating from the corresponding governing equations. This \jpang{idea} is the core of feedback optimization.
			Then, the result obtained by the feedback optimization algorithm is sent back to the physical system to regulate its power generation. Consequently, this process constitutes a feedback loop, leading to the concept of feedback optimization.}
		
		\begin{figure}[t]
			\centering
			\includegraphics[width=0.38\textwidth]{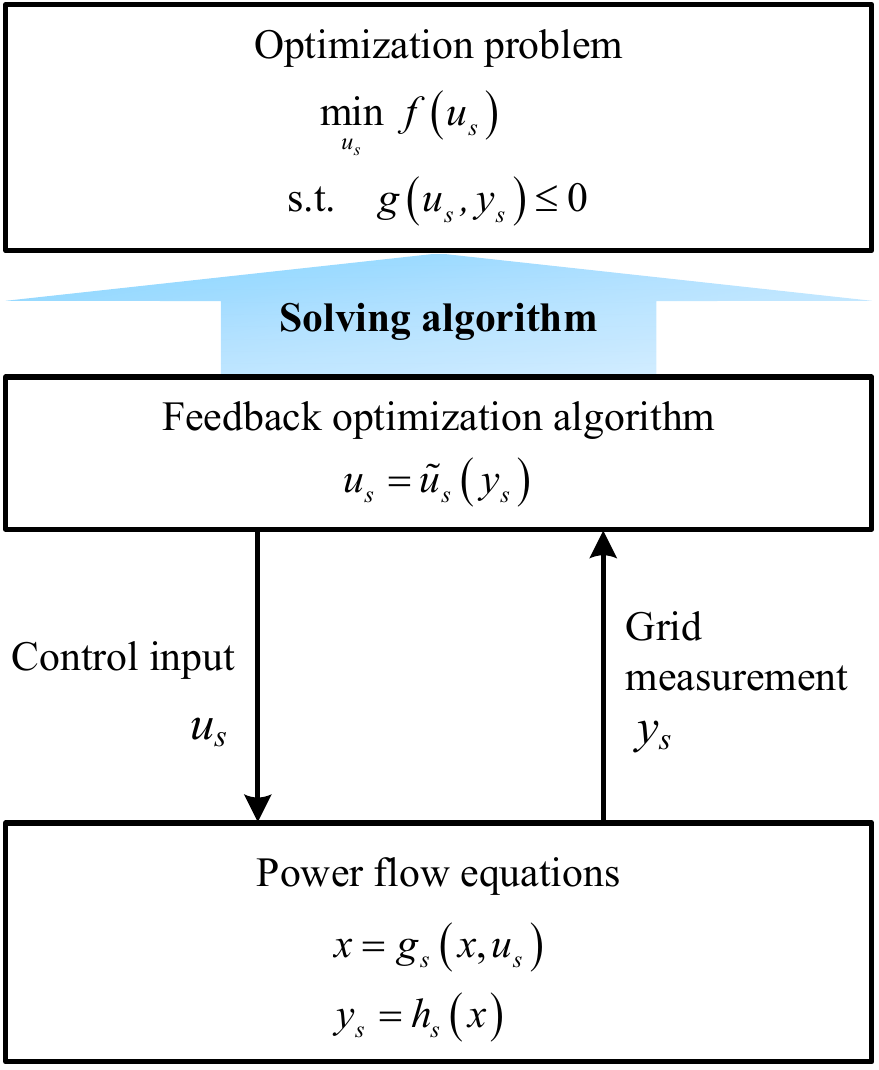}
			\caption{Diagram of feedback optimization}
			\label{Feedback_optimization}
		\end{figure}
		
		\jpang{Fig.\ref{Feedback_optimization} presents the conceptual diagram of feedback optimization. The lower level is the (quasi)steady-state system model described by the power flow equations, which is introduced in Section \ref{Power_flow_equations}. Here we give its simplified notation. 
			\begin{align*}
				0 = {g_s}\left( {{ x},{{ u}_s}} \right)
			\end{align*} 
			where $  x $ is the system state, ${ u}_s$ is the controllable variable, and ${ y}_s$ is the system output. For example, ${ u}_s$ could be the generation set point and ${ y}_s$ could be the voltage and power generation. If we consider the measurement $y_s$ as the parameter, the single-period problem has the general form. 
			\begin{subequations}\label{Feedback_algorithm_problem}
				\begin{align}
					\mathop {\min }\limits_{{{ u}_s}}\quad & f\left( {{{ u}_s}} \right)\\
					\text{s.t.}\quad &g\left( {{{ u}_s},{{ y}_s}} \right) \le 0\\
					&\text{Power flow equations}
				\end{align}
			\end{subequations}
			The power flow equations serve as a set of equality constraints, which, however, do not need to be solved mathematically. By measuring ${y}_s$ from the physical power system, this constraint is satisfied naturally, and can therefore be removed. Thereafter, some iterative algorithm is designed to solve \eqref{Feedback_algorithm_problem} together based on the measurement ${y}_s$. The intermediate iterates are then sent back to the physical system to drive it to the desired working point. The new measurement will be obtained to continue to implement the next iteration.}
		
		\jpang{Feedback optimization has two salient advantages. 1) Low computational complexity: the measurement remarkably reduces the computational complexity of solving the high-dimensional nonlinear power flow equations. 2) Fast response speed: feedback optimization has a fast response speed and is adaptable to the time-varying environment caused by volatile load demand and renewable generation sources. }
		
		\jpang{Here, we take the online Optimal Power Flow (OPF) problem in \cite{gan2016online} as a concrete example and illustrate how to design feedback optimization algorithms.}
		
		\begin{example}
			OPF is a fundamental problem in power system operations, which has broad applications \cite{huneault1991survey}, including economic dispatch, demand response, state estimation, stability assessment, etc. Substantial efforts have been devoted to finding reliable solution methods to OPF since it was first formulated in 1962 by Carpentier \cite{carpentier1962contribution}. 
			The OPF problem is formulated as the following optimization problem: 
			\begin{subequations}\label{OPF_Problem}
				\begin{align}\label{OPF_Objective}
					\min &\quad f(x)=\sum_{i=0}^{n} a_{i} P_{i}^{2}+b_{i} P_{i}\\
					\text { over }& x :=\left(P_{i}, Q_{i}, i \in \mathcal N\right)\nonumber \\
					&	y :=\left(P_{0}, Q_{0}, v_{i}, i \in \mathcal N ; P_{i j}, Q_{i j}, \ell_{i j},(i, j) \in\mathcal E\right)\nonumber
					\\
					s.t.\ &\eqref{DistFlow}\\
					&\label{voltage_OPF}
					\underline{v}_{i} \leq v_{i} \leq \overline{v}_{i}, \quad i \in \mathcal N\\
					&\label{P_limits}
					\underline{P}_{i} \leq P_{i} \leq \overline{P}_{i} , \quad  i \in \mathcal N\\
					&\label{Q_limits}
					\underline{Q}_{i} \leq Q_{i} \leq \overline{Q}_{i}, \quad  i \in \mathcal N
				\end{align}
			\end{subequations}
			where $ a_{i}>0, b_i>0 $ are constants. The objective function is to minimize the generation cost, with constraints \eqref{DistFlow} and \eqref{voltage_OPF} representing the power flow equation and voltage limits respectively, and \eqref{P_limits} and \eqref{Q_limits} representing the power injection constraints. 
			
			Define the domain set
			\begin{align*}
				X:=\left\{\left(P_{i}, Q_{i}\right) \mid \underline{P}_{i} \leq P_{i} \leq \overline{P}_{i}, \underline{Q}_{i} \leq Q_{i} \leq \overline{Q}_{i}, i \in \mathcal N\right\}
			\end{align*}
			Then, we can define a continuously differentiable function $F$ to rewrite the power flow equation \eqref{DistFlow} as follows.
			\begin{align}
				F(x, y)=0
			\end{align}
			In fact, $ F(x, y)=0 $ defines implicitly a function $y=y(x)$ over $X$ :
			\begin{align*}
				P_{0}&:=P_{0}(x), Q_{0}:=Q_{0}(x) & \\
				v_{i}&:=v_{i}(x), \qquad\qquad i \in \mathcal N \\
				P_{i j}&:=P_{i j}(x), Q_{i j}:=Q_{i j}(x), \ell_{i j}=\ell_{i j}(x), \quad (i,j)\in \mathcal E
			\end{align*}
			Thus, the OPF problem \eqref{OPF_Problem} can be written equivalently as
			\begin{subequations}\label{OPF_Problem2}
				\begin{align}
					\mathop {\min }\limits_{x\in X} \quad & a_{0} P_{0}^{2}(x)+b_{0} P_{0}(x)+\sum_{i=1}^{n}\left(a_{i} P_{i}^{2}+b_{i} P_{i}\right)\\
					s.t.\quad& \label{v_constraints} \underline{v}_{i} \leq v_{i}(x) \leq \bar{v}_{i}, \quad i \in \mathcal N
				\end{align}
			\end{subequations}
			The equivalence between \eqref{OPF_Problem} and \eqref{OPF_Problem2} is based on the implicit assumption that $ F(x, y)=0 $ always holds. This is further explained in the online algorithm that follows.
			
			In \cite{gan2016online}, to simplify the feasible set, the constraints \eqref{v_constraints} is added to the objective function via a log-barrier function. Then, the Lagrangian of \eqref{OPF_Problem2} is 
			\begin{align}
				&L(x, \underline\mu, \bar \mu):= a_{0} P_{0}^{2}(x)+b_{0} P_{0}(x)+\sum_{i=1}^{n}\left(a_{i} P_{i}^{2}+b_{i} P_{i}\right)\nonumber \\
				&\qquad-\underline{\lambda} \sum_{i=1}^{n} \ln \left(v_{i}(x)-\underline{v}_{i}\right)-\bar{\lambda} \sum_{i=1}^{n} \ln \left(\bar{v}_{i}-v_{i}(x)\right)
			\end{align}
			where $ \underline{\lambda}, \overline{\lambda} $ are decreasing sequences such that 
			\begin{align*}
				&\lim \limits_{v \rightarrow \underline{v}^+_{i}}-\underline{\lambda} \ln \left(v-\underline{v}_{i}\right)=\infty, \ i \in \mathcal N\\
				& \lim \limits_{v \rightarrow \overline{v}^-_{i}}-\overline{\lambda} \ln \left(\overline{v}_{i}-v\right)=\infty, \ i \in \mathcal N
			\end{align*}
			

			The online OPF algorithm is proposed based on the gradient projection method.
			\begin{subequations}\label{Online_OPF_Algorithm}
				\begin{align}\label{Online_OPF_Iteration}
					x_{t+1} &=\mathcal{P}_{{X}}\left[x_{t}-\eta \nabla_{x} L(x, \underline\lambda, \bar \lambda)\right] \\
					\label{Online_OPF_Measurement}
					y_{t} &=y(x_{t})
				\end{align}
			\end{subequations}
			where $\eta>0$ is the stepsize, and $ \mathcal{P}_{{X}}[\cdot] $ is the projection onto $X$. The gradient equals
			\begin{align*}
				&\nabla_{x} L(x, \underline\lambda, \bar \lambda) =\partial_{x} L(x, y(x) , \underline\lambda, \bar \lambda)\\
				&\qquad-\partial_{y} L(x, y(x) , \underline\lambda, \bar \lambda){\left[\partial_{y} F(x, y(x))\right]^{-1} \partial_{x} F(x, y(x)) }
			\end{align*}
			where $ \partial_{x}L $ stands for the partial derivative of $ L $ with respect to $ x $. It is shown in in \cite{gan2016online} that $ \nabla_{x} L(x, \underline\lambda, \bar \lambda) $ can be approximated by the linearized BFM \eqref{Linear_DistFlow}, which reduces the computational burden greatly.

			\begin{algorithm}[t]
				\caption{Flow chart of feedback-based OPF }
				\label{OPF_algorithm}
				{\color{black}
					\begin{algorithmic}[1]
						\State Initialization: Set $t=1$, the value of $ x(1) $, parameter $ \eta$, and decreasing sequences $\underline\lambda, \bar\lambda $.
						\While{$t\le T$}
						\State Apply the current iterate $x_{t}$ to the physical system.
						\State Measure $y_{t}$ directly from the network.
						\State Compute $x_{t+1}$ from \eqref{Online_OPF_Iteration} and objective function from \eqref{OPF_Objective}.
						\State Update $ t\leftarrow t+1 $.
						\EndWhile
						\State Output: $ x_t, f_t(x_t), t=1, \cdots, T $.
				\end{algorithmic}}
			\end{algorithm}
			
			\jpang{The online implementation of \eqref{Online_OPF_Algorithm} is introduced and represented in Algorithm \ref{OPF_algorithm}.
				In step 2, the power flow equation is solved by the physical system based on the law of physics. In this way, we avoid solving the power flow equation \eqref{DistFlow}, which is computationally costly. The algorithm could continue with the much simpler measurement and feedback steps. }
			
		\end{example}
		
		
		\subsection{Application survey}\label{Feedback_optimization_application}

		
		In this subsection, we will introduce more applications of feedback optimization, starting with the online OPF problem.
		Works on online OPF could be roughly divided into three types based on the methods adopted, including gradient descent: \cite{gan2016online,dall2018optimal,ortmann2020experimental,picallo2021adaptive}, alternating direction method of multipliers (ADMM)  \cite{zhang2017dynamic}, and the Newton-Raphson method \cite{tang2017distributedCDC,tang2017real,picallo2021cross}. 
		
		In the first type, the primal-dual gradient method is used to solve the time-varying OPF problems \cite{dall2018optimal,ortmann2020experimental}. The gradient approximation method could be adopted to reduce its computation burden \cite{gan2016online}. {\color{black}If some buses are not measurable, sensitivity or Kalman filter could be utilized to estimate variables in these buses, as was done in \cite{picallo2020closing,picallo2021adaptive}. }
		Dynamic ADMM is similarly adopted in \cite{zhang2017dynamic} to solve the OPF problem based on online measurements. 
		To accelerate the convergence speed, Newton-Raphson methods are designed  \cite{tang2017distributedCDC,tang2017real,picallo2021cross}, which utilizes the second-order information. The second-order information helps to better handle the non-linearity of the OPF problem, so these methods usually have improved tracking performance.

		
		Another important application is online voltage control.
		Generally speaking, optimal voltage regulation can be viewed as a type of OPF problem, which covers different kinds of objectives, such as minimizing power losses and voltage difference from nominal value \cite{zhu2016fast, liu2017decentralized, zhou2018reverse}. 
		{\color{black}Due to the multi-objective feature, typically weighted sums of these cost functions are considered. The branch model for the distribution system, i.e., the Linearized DistFlow in \eqref{Linear_DistFlow}, is widely utilized in the problem formulation. }
		Recent works can be roughly categorized into two types: decentralized control and distributed control. For the decentralized voltage control, locally available information such as bus voltage magnitude is utilized to design the controller \cite{turitsyn2011options}. The control objective associated with the voltage error has to be defined in the weighted form, i.e., the admittance matrix-induced norm \cite{zhu2016fast, liu2017decentralized}. In this way, its gradient equals the nodal power flow equation, which spawns the online measurement and local implementation. The specific form of the objective function also restricts the extensibility of the local control. 
		The distributed voltage control can avoid the disadvantage to some extent, while the communication is needed for the tradeoff \cite{antoniadou2017distributed}, such as communication with immediate neighbors \cite{bolognani2013distributed, vsulc2014optimal, bolognani2015distributed, zhang2015optimal, liu2018distributed, liu2018hybrid} or two-hop neighbors \cite{Tang2018Fast,wang2020asynchronous_voltage}. Then, the objective could be more general, such as the voltage error in the Euclidean norm and the total power losses.

		\begin{remark}[Comparison between OGDC and feedback optimization]
			Both OGDC and feedback optimization are designed for the power system operation in the time-varying environment. The major similarity between OGDC and feedback optimization is that they both adopt the idea of feedback control. However, they also have the following differences in the following two aspects.
			\begin{itemize}
				\item Feedback optimization focuses on optimization at a slower time scale with the consideration of (quasi)steady states, while the OGDC considers the dynamics of physical systems on a faster time scale. 
				\item Feedback optimization is more concerned with operation optimality, where tracking errors compared with offline results are taken as the performance index. OGDC pays more attention to the system stability, where optimality of the control scheme refers to the result in the steady state. 
			\end{itemize}  
			Nevertheless, in practical operation, one sometimes needs to combine two methods to achieve satisfactory performance on different time scales. 
		\end{remark}

		\section{Online Optimization for Multi-Period Problems}\label{sec:mpo}
		In this section, we introduce two commonly-used online optimization methods for multi-period problems, including Lyapunov optimization and OCO. We will present their main ideas and typical applications. 
		Here, we emphasize a pure online decision-making process without referring to historical data or prediction. 
		
		\subsection{Lyapunov Optimization}\label{Lyapunov_Optimization}
		In this subsection, we first introduce the main idea of Lyapunov optimization for multi-period problems. Then, we survey its applications in power systems. 
		
		\subsubsection{Main idea} \label{Lyapunov_main_idea}
		In a multi-period optimization problem with $ T $ periods, there exists a queue defined as 
		\begin{align}\label{Def_queue}
			Q_{t+1}=Q_{t}+x_t,\quad t \in\{0, \ldots, T-1\}
		\end{align}
		where $ Q_t $ is the state variable. $ x_t=Y_1(w_t,u_t) $ is assumed to be bounded, where $ w_t $ is a stochastic variable, $ u_t $ is the decision variable, $ Y_1(\cdot) $ is a mapping from input to state determined by the system property. In power systems, $ Q_t $ is usually the SoC of an ESS, and then $ x_t $ is the charging or discharging rate of the ESS. 
		
		Suppose the initial value $Q(0)$ is a constant. Summing up both sides of \eqref{Def_queue} over $t=0, \cdots, T-1$ and taking expectation give rise to
		\begin{align}\label{Sum_queue}
			\mathbb{E}\left[Q_T\right]-Q(0)=\sum\nolimits_{t=0}^{T-1} \mathbb{E}\left[x_t\right]
		\end{align}
		Assuming $Q_T$ is bounded, the left-hand side of \eqref{Sum_queue} is also bounded. Then, dividing both sides of \eqref{Sum_queue} by $T$ and taking limits as $T \rightarrow \infty$ yields
		\begin{align}\label{Mean_rate_stability}
			\lim _{T \rightarrow \infty} \frac{1}{T} \sum_{t=0}^{T-1} \mathbb{E}\left[x_t\right]=0
		\end{align}
		In \eqref{Mean_rate_stability}, it shows that net storage in $Q$ is zero over a long-term horizon, which is the so-called mean-rate stability.
		
		We focus on the following multi-period problem
		\begin{subequations}\label{Optimization_problem}
			\begin{align}
				\min _{u_t,\forall t}\quad&  \lim _{T \rightarrow \infty} \frac{1}{T} \sum_{t=0}^{T-1} \mathbb{E}[f(y_t)]\\
				\text{s.t.}\quad&  g(y_t)\le0\\
				& h(y_t)=0\\
				& \text{the queue definition}\ \eqref{Def_queue}\nonumber\\
				& \text{the mean rate stability}\ \eqref{Mean_rate_stability}\nonumber
			\end{align}
		\end{subequations}
		where $ y_t=Y_2(w_t,u_t) $ with a mapping $ Y_2(\cdot) $, which is not necessarily equal to $ x_t $. 
		Problem \eqref{Optimization_problem} is to minimize the long-term time-average cost with equality and inequality constraints. Since there exists future variables $ x_{t+1}, x_{t+2}, \cdots $ at any time slot $t$, \eqref{Optimization_problem} is difficult to be solved online. 
		The philosophy of Lyapunov optimization requires making an online decision only based on $x_t$ and $w_t$ observed in the current stage, as well as the current queue $Q_t$. Thus, the key is to eliminate \jpangk{dependency on the} future variables in the problem formulation.
		
		We first define a quadratic Lyapunov function
		\begin{align}
			L_{t}=\frac{1}{2} Q_t ^{2}
		\end{align}
		The Lyapunov drift between two adjacent time slots is
		\begin{align}
			\Delta_{t}=L_{t+1}-L_{t}=\frac{Q^{2}_{t+1}-Q^{2}_t}{2}
		\end{align}
		
		It is verified in \cite{neely2010stochastic} that minimizing the Lyapunov drift ensures the stability of $Q_t$. If we further  consider minimizing the cost together, the objective function could be drift-plus-cost, i.e., $\Delta_{t}+\rho f(x_t)$, where $\rho$ is a constant adjusting the weight between drift and cost. Because $Q_{t+1}$ in $\Delta_{t}$ is not known in the time slot $t$, the drift-plus-cost function cannot be minimized directly. Recalling the definition in \eqref{Def_queue}, we have $ Q^2_{t+1}= \left(Q_t+x_t\right)^2 $. Then, we get the upper bound of $\Delta_{t}$.
		\begin{align}
			\Delta_{t}&=\frac{x^2_t}{2}+Q_tx_t\le B+Q_tx_t
		\end{align}
		where $ B=\frac{\overline{x}^2}{2} $ is a constant with $ {\overline{x}} $ as the upper bound of $ x_t, \forall t $. The upper bound is only determined by the information of the current stage. Taking $ Q_t\cdot x_t+\rho\cdot f(y_t) $ as the objective function and removing the mean rate stability constraint, the original problem \eqref{Optimization_problem} is relaxed to 
		\begin{subequations}\label{Optimization_problem_drift}
			\begin{align}\label{Optimization_problem_drift_obj}
				\min _{u_t,\forall t}\quad&  Q_t\cdot x_t+\rho\cdot f(y_t)\\
				\text{s.t.}\quad&  g(y_t)\le0\\
				& h(y_t)=0
			\end{align}
		\end{subequations}
		{\color{black}Since $ B $ is a constant, removing it in \eqref{Optimization_problem_drift_obj} does not change the solution of the problem. Equation \eqref{Def_queue} is satisfied by the physical law, and thus is also not included explicitly.} Clearly, all the parameters in problem \eqref{Optimization_problem_drift} can be obtained at time slot $t$, which can now be solved directly.
		
		The main steps of the Lyapunov optimization to solve \eqref{Optimization_problem} up till now is summarized in Algorithm \ref{Lyapunov_algorithm}.
		
		\begin{algorithm}[t]
			\caption{Online Algorithm based on Lyapunov Optimziation }
			\label{Lyapunov_algorithm}
			{\color{black}
				\begin{algorithmic}[1]
					\State {Initiation}: Set $t=1$, the value of $ Q_1 $ and parameter $ \rho $. 
					\State Transform \eqref{Optimization_problem} to the online counterpart \eqref{Optimization_problem_drift}.
					\While{$t\le T$}
					\State S1: Observe the stochastic variable $w_t$.
					\State S2: Solve problem \eqref{Optimization_problem_drift} to get $u_t$ and objective function.
					\State S3: Get $x_t$ and $y_t$ from the mapping $ Y_1(\cdot), Y_2(\cdot) $, respectively.
					\State S4: Update the queue $ Q_{t+1} $ according to \eqref{Def_queue}.
					\State S5: Update $ t\leftarrow t+1 $.
					\EndWhile
					\State Output: $ x_t, y_t, f_t(y_t), Q_t, t=1, \cdots, T $.
				\end{algorithmic}
			}
		\end{algorithm}
		
		One can show that Algorithm \ref{Lyapunov_algorithm} solves the relaxed problem without the constraint of mean-rate stability. Further, it is proved in \cite[Chapter 4.1]{neely2010stochastic} that the queue obtained from Algorithm \ref{Lyapunov_algorithm} is mean-rate stable. Moreover, the error between the result obtained by Algorithm \ref{Lyapunov_algorithm} and the optimal solution of \eqref{Optimization_problem} is bounded, which is determined by the parameter $\rho$. Generally, a larger $\rho$ leads to a smaller error \cite[Chapter 4.1]{neely2010stochastic}. 
		
		The Lyapunov optimization also can be interpreted as the stochastic dual gradient method \cite{Tianyi2017Averaging}, where the queue in \eqref{Def_queue} is the iteration of the Lagrange multiplier. From this perspective, some approaches could be adopted to improve the performance of the solution obtained by Lyapunov optimization. For example, an extra gradient evaluation is added in \cite{chen2018learn}, which learns from the historical data and then adapts to the upcoming strategies. 
		Because the procedure of Algorithm \ref{Lyapunov_algorithm} is easy to follow, we do not introduce a specific example here.

		
		%

		\subsubsection{Application survey}
		
		The operation of an ESS naturally fits the queue given in \eqref{Def_queue}, where $ Q_t $ is the SoC, and $ x_t $ is the sum of charging and discharging power. Thus, the online algorithm based on the Lyapunov optimization has been mainly reported in storage-related scenarios \cite{stai2021online}. A typical application is the online economic dispatch, including energy management in smart grid with distributed energy resources  \cite{salinas2013dynamic,fan2020online}, in \cite{huang2014adaptive,Wenbo2015real,paul2020real,hao2020distributed,zeinal2020flexible,shotorbani2021enhanced} for microgrids. Moreover, it can also be utilized in integrated energy systems \cite{zhou2014optimal,zhang2018online,li2020lyapunov,zeinal2020real,wang2021event}, where the queue could be the SoC of heat storage or the room temperature. A similar method is also applied in \cite{yan2019real,ahmad2019real} for the charging of electric vehicles, where the SoC of the EV battery is treated as the queue. Besides, another important application lies in power markets, such as the energy sharing with storage systems \cite{liu2017online,liunian2017online,dafeng2020energy}, where each participant has its own ESS. 
		
		The aforementioned works focus on the centralized Lyapunov optimization, which, however, also can be implemented in the distributed manner \cite{zheng2014distributed,sun2016distributed,zhong2019online,fan2020optimal,zhu2021distributed}. In \cite{zheng2014distributed}, the heating ventilation, and air-conditioning system are considered in DR, where the controlled room temperature is modeled as a queue similar to the SoC of a battery. The algorithm is implemented in a partially distributed way with a control center broadcasting the summary of power consumption. In \cite{sun2016distributed,zhong2019online}, the alternating direction method of multipliers (ADMM) is adopted in the distributed implementation, which is also a partially distributed way with an aggregator updating dual variables. In \cite{fan2020optimal}, the dual ascent algorithm is designed to coordinate networked DERs as a virtual power plant in a distributed manner, which is also partially distributed. In these works, a central coordinator is needed to broadcast control parameters or queue states \cite{zhu2021distributed}, and a fully distributed manner still needs to be developed. 
		
		In most works, the time horizon $ T $ is assumed to be infinite, which can simplify the proof of the algorithm performance mathematically. Nonetheless, we only focus on the energy management during a given period in many cases, which requires the finite-time horizon Lyapunov optimization \cite{li2015real,li2016real,li2018residential}. 
		In \cite{li2015real}, a real-time is designed for the management of batteries within a finite period. A similar method is also utilized in \cite{li2016real,li2018residential}, considering joint energy storage management and load scheduling. The procedure of finite-time Lyapunov optimization is similar to the infinite counterpart, which makes the performance proof more difficult. 
		
		\jpang{Lyapunov optimization is an effective tool for multi-period decision-making problems in power systems with time-coupled states, such as the SoC of the energy storage systems, or the room temperature of HVAC or district cooling systems. Because it performs an ``1-lookahead" mechanism, only a single-period deterministic problem needs to be solved at each period. As such, an analytical solution to the Lyapunov-drift problem \eqref{Optimization_problem_drift} becomes possible. For example, in \cite{Zhongjie2021Real}, a multiparametric programming method is used to get the analytical real-time dispatch policy. }

		\subsection{Online Convex Optimization}\label{Online_Convex_Optimization}
		
		From Section \ref{Lyapunov_Optimization} and Algorithm \ref{Lyapunov_algorithm}, the stochastic variable $w_t$ needs to be observed first before solving \eqref{Optimization_problem_drift}. This is the so-called ``1-lookahead" pattern, which covers some problems in power systems. For example, the electricity price is sent to the wind-storage integrated power plant, after which the self-dispatch action is taken \cite{Zhongjie2021Real}. However, this is not always the case. In some cases, the uncertainties are unknown before the decision-making and will be revealed later, which leads to the ``0-lookahead" pattern. For example, a market participant, such as the ESS owner, bids on a capacity without the knowledge of the electricity price. The price could be obtained only after the market is cleared. In this circumstance, Lyapunov optimization is inapplicable, and instead, OCO is typically used. In this section, we will introduce its main idea and applications in power systems.
		
		\subsubsection{Main idea}
		OCO focuses on the following problem.
		\begin{subequations}\label{OCO_round_problem}
			\begin{align}
				\min _{x_{t}}&\quad  f_{t}\left(x_{t}\right) \\
				\text { s.t. }&\quad g_{t}\left(x_{t}\right) \le0, \quad t=1, \ldots, T
			\end{align}
		\end{subequations}
		where the convex objective function $f_{t}(\cdot)$ and the convex constraint $g_{t}(\cdot)$ are unknown before $ x_{t} $ is determined by OCO algorithm, i.e., uncertainties not observed in prior. It is usually viewed as a repeated game between a learner and nature across a finite time horizon $t=1, \ldots, T$, which may be adversarial \cite{zinkevich2003online}. At the beginning of each time slot $t$, the learner determines a decision variable $x_{t}\in X$ by a pre-designed online algorithm. Then, nature selects the value of uncertainty, and thus determines $f_{t}(\cdot)$ and $g_{t}(\cdot)$. Consequently, the learner gets the objective $f_{t}\left(x_{t}\right)$ and the constraint $g_{t}\left(x_{t}\right) \leq 0$. 
		
		Most of the pre-designed algorithms are online variants of the discrete form of primal-dual gradient method. 
		\begin{subequations}\label{Primal_dual}
			\begin{align}\label{Primal_update}
				x_{t} &=\mathcal{P}_{{X}}\left[x_{t-1}-\eta \left(\nabla_{x} f_{t-1}\left(x_{t-1}\right)+\lambda_{t-1}^T\nabla_{x} g_{t-1}\left(x_{t-1}\right)\right)\right] \\
				\label{Dual_update}
				\lambda_{t} &=\mathcal{P}_{\mathbb{R}^+}\left[\lambda_{t-1}+\eta g_{t-1}\left(x_{t-1}\right)\right]
			\end{align}
		\end{subequations}
		where $ \lambda_{t} $ is the Lagrangian multiplier with respect to the constraint $g_{t}\left(x_{t}\right) \leq 0$, $\eta$ is a constant stepsize, and $ {\mathbb{R}^+} $ is the nonnegative Euclidean space with proper dimension. The OCO algorithm based on \eqref{Primal_dual} is given in Algorithm \ref{OCO_algorithm}, which clearly shows the ``1-lookahead" pattern. 
		
		\begin{algorithm}[t]
			\caption{Online Algorithm based on Lyapunov Optimziation }
			\label{OCO_algorithm}
			{\color{black}
				\begin{algorithmic}[1]
					\State {Initiation}: Set $t=1$, the value of $ x_0, \lambda_{0} $ and parameter $ \eta $. \eqref{Optimization_problem_drift}.
					\While{$t\le T$}
					\State S1: Update the primal variable $ x_t $ by \eqref{Primal_update}.
					\State S2: Update the dual variable $ \lambda_t $ by \eqref{Dual_update}.
					\State S3: Determine the corresponding objective function $ f_t(x_t) $ and constraint $ g_t(x_t) $.
					\State S4: Compute the performance index, $ \operatorname{Reg}$, $\operatorname{Vio}$, $\operatorname{CR} $.
					\State S5: Update $ t\leftarrow t+1 $.
					\EndWhile
					\State Output: $ x_t, f_t(x_t), t=1, \cdots, T $, $ \operatorname{Reg}(T)$, $\operatorname{Vio}(T)$.
				\end{algorithmic}
			}
		\end{algorithm}

		In \eqref{Primal_dual}, because $x_{t}$ is obtained with the information of  $f_{t-1}(\cdot), g_{t-1}(\cdot)$ instead of $f_{t}(\cdot), g_{t}(\cdot)$, it is almost impossible to obtain the optimal solution to \eqref{OCO_round_problem}.
		A natural question is how to assess the performance of the pre-designed algorithm, i.e., the suboptimality of $x_{t}$. Here, we introduce three commonly used measures: Regret, Violation, and Competitive Ratio. 
		
		The first measure \emph{regret} is defined as
		\begin{align*}
			\operatorname{Reg}^d(T):=\sum_{t=1}^{T}\left[f_{t}\left(x_{t}\right)-f_{t}\left(x_{t}^{*}\right)\right]
		\end{align*}
		where $ x_t^* $ is the optimal solution to the problem \eqref{OCO_problem}. Regret quantifies the suboptimality caused by the pre-designed algorithm. Since $ x_t^* $ is the optimum at time slot $t$, $ \operatorname{Reg}^d(T) $ is sometimes called \emph{dynamic regret}. In contrast, we can also define the \emph{static regret}.
		\begin{align*}
			\operatorname{Reg}^s(T):=\sum_{t=1}^{T}\left[f_{t}\left(x_{t}\right)-f_{t}\left(x_s^{*}\right)\right]
		\end{align*}
		where $ x_s^{*} $ is the optimal solution to \eqref{OCO_problem}.
		\begin{subequations}\label{OCO_problem}
			\begin{align}
				\min _{x}&\ \sum_{t=1}^{T} f_{t}\left(x_{t}\right) \\
				\text { s.t. }&\ g_{t}\left(x_{t}\right) \leq 0,\ t=1, 2, \cdots, T
			\end{align}
		\end{subequations}
		The baseline strategy $ x_s^* $ in the static regret remains identical throughout the horizon, which may cause large deviations between $x_t$ and $ x_s^* $ with the increase of $ t $ in the time-varying environment. In contrast, the baseline strategy $ x_t^* $ in the dynamic regret varies with the time slot $t$, which seems to be more suitable for performance guarantee in power systems with strong volatility. 
		
		{\color{black}The second measure \emph{violation} is defined as  
			\begin{align}
				\operatorname{Vio}(T):=\sum_{t=1}^{T}\left\|\left[ g_{t}\left(x_{t}\right)\right]_{+}\right\|
			\end{align}
			where the operator $ [a]_+:=\max(a,0) $. Solutions that are feasible to the set of constraints do not contribute to $ \operatorname{Vio} $. 
			
			Finally, the last measure is the \emph{Competitive Ratio} $ \operatorname{CR} $, which is the rate between objective functions obtained by the online algorithm and optimal value \cite{Yanfang2021optimal,menati2022competitive}. }
		\begin{align}
			\operatorname{CR}=\max _{\forall t} \frac{f_{t}\left(x_{t}\right)}{f_{t}\left(x_{t}^{*}\right)}
		\end{align}
		Clearly we have $\operatorname{C R}\geq 1$. This index is to drive the cost of the online algorithm close to the offline optimum.

		
		In existing works, many extensions and variants of the OCO algorithm are investigated, which includes updating the primal and dual variables in a Gauss-Seidel manner \cite{Tianyi2017An}, adding a regularization term to get the strong concavity of the dual variable \cite{cao2018online}, distributed implementation of the algorithm \cite{Xinlei2020Distributed,lesage2020dynamic,Deming2022Distributed}, etc. 
		In addition to the iteration type, another important variant of the OCO algorithm is based on the information obtained by the learner. If the analytical expressions of $f_{t}(\cdot)$ and $g_{t}(\cdot)$ or their gradients are available, this is called full feedback. Otherwise, it is known as the bandit feedback case, i.e., the values of $f_{t}(\cdot)$ and $g_{t}(\cdot)$ obtained only at the sampling instance \cite{flaxman2005online,cao2018online,Xinlei2021Distributed}. In the case of bandit feedback, the key is to estimate the gradient of $ f_t(x_t), g_t(x_t) $ using limited information. 
		It must be pointed out that the OCO algorithm is strongly problem-dependent, and has no unified mathematic paradigm. Comparisons between the different OCO algorithms can be found in \cite{hazan2016introduction}. 
		Since the procedure of Algorithm \ref{OCO_algorithm} is easy to follow, we do not introduce a specific example here. 

		\subsubsection{Application survey}
		
		The framework of OCO is first defined in the machine learning literature \cite{shalev2011online,hazan2016introduction}, which has recently gained attention in power systems, particularly in the applications of demand-side management \cite{badiei2015online,ma2016distributed,kim2016online,bahrami2017online,bahrami2019online,lesage2019online,lesage2020dynamic,chen2021online}. 
		In \cite{badiei2015online}, an online algorithm is studied to address optimization problems with ramp constraints, which presents asymptotically tight bounds on the competitive difference. 
		In \cite{kim2016online}, the varying price elasticity of consumers in the DR is considered by the online algorithms, where both full and bandit feedback structures are included. 
		In \cite{bahrami2017online}, the long-term load scheduling problem is investigated, which is modeled as a partially observable stochastic game due to the uncertainties of price and load demand. To solve the problem and get the Markov perfect equilibrium, an online load scheduling learning algorithm is proposed based on the actor-critic method. 
		Similar methods are extended in  \cite{bahrami2019online} to deploy the DR programs for data centers. 
		In \cite{chen2021online}, the residential DR is formulated as a contextual multi-armed bandit problem, which is solved by an online learning and selection algorithm based on the Thompson sampling method.

		The OCO algorithm also can be implemented in a distributed manner. In \cite{ma2016distributed}, distributed online learning algorithm is proposed for the charging control of electric vehicles, which requires only one-way communication, i.e., the distribution company publishes the prices. This structure fits in with the current communication infrastructure in reality.
		In \cite{lesage2020dynamic}, the DR considering heating, ventilation, and air-conditioning systems of commercial buildings is investigated, where the dynamic regret is used to evaluate the performance of the online distributed weighted dual averaging algorithm.
		Frequency regulation using ESS based on the OCO is studied in \cite{zhao2020distributed}, which also can be realized in a distributed manner with the capability of plug-and-play.
		
		\jpang{In an OCO framework, the objective functions and constraints are revealed after the decisions are made, which contradicts most cases in traditional power systems. Consequently, the application of OCO in power systems is still limited currently. Nonetheless, with the high penetration of volatile renewable generations, the OCO method shows great potential in the future due to its features, such as its fast response speed and its independence from the need for predictions. }
		
		\begin{remark}[Comparison between Lyapunov optimization and OCO]
			Lyapunov optimization and OCO are both widely used in multi-period problems. Their differences are twofold.
			\begin{itemize}
				\item The major difference is the decision-making process. Lyapunov optimization is the ``1-lookahead" pattern, where the uncertainty $w_t$ is observed first, and then the Lyapunov drift problem \eqref{Optimization_problem_drift} will be solved. In contrast, OCO is usually a ``0-lookahead" pattern, where the decision is made before the observation of the uncertainty. This difference is illustrated in Fig.\ref{Lyapunov_OCO}.
				\begin{figure}[t]
					\centering
					\includegraphics[width=0.48\textwidth]{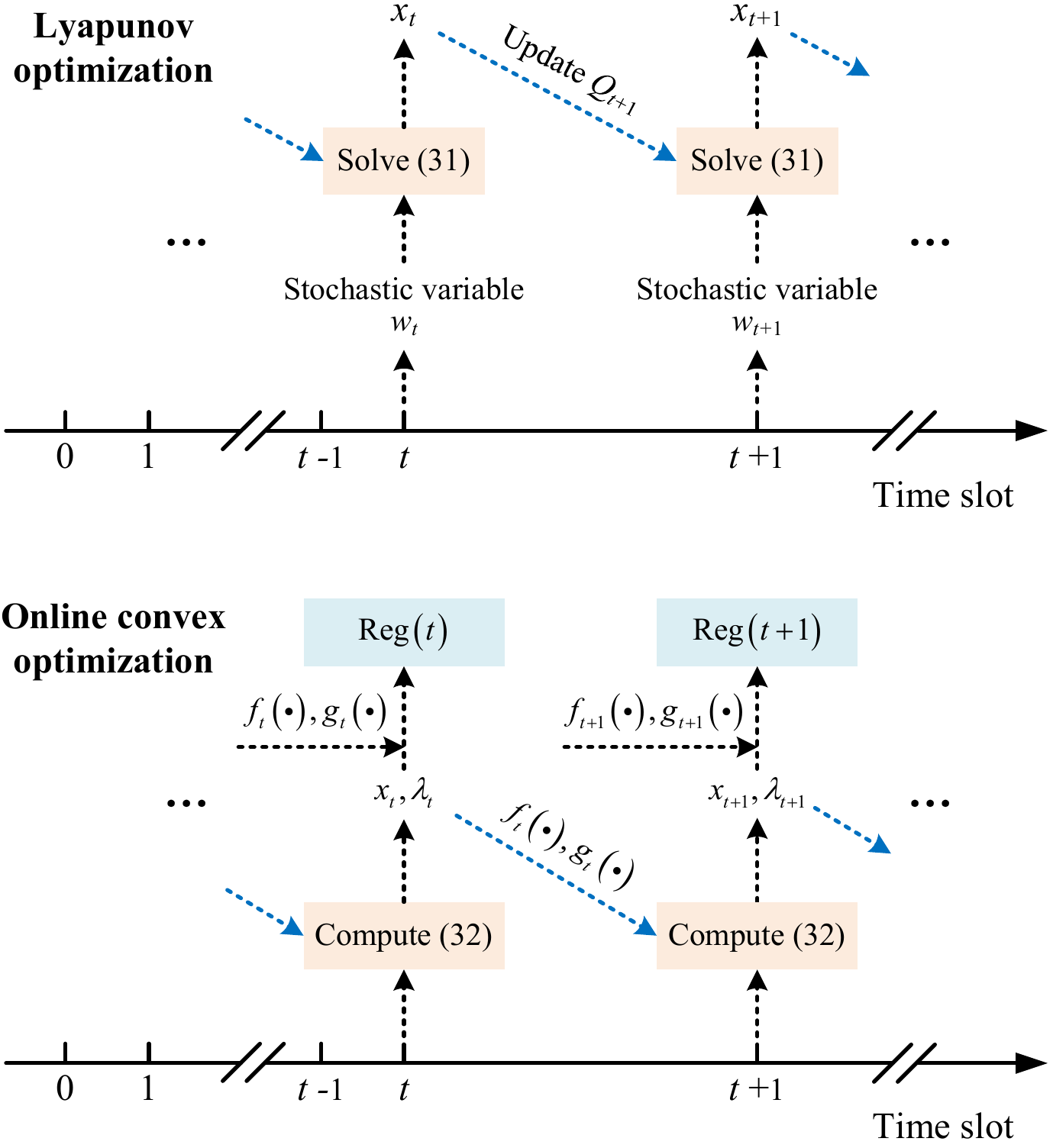}
					\caption{Decision-making process of Lyapunov optimization and
						online convex optimization}
					\label{Lyapunov_OCO}
				\end{figure}
				\item Lyapunov optimization usually follows a standardized paradigm in the algorithm design, with attention paid to the problem formulation. On the contrary, OCO algorithms are strongly problem-dependent without a unified mathematical paradigm, which, however, all focus on the performance estimation with metrics like $ \operatorname{Reg}$, $\operatorname{Vio}$, and $\operatorname{CR} $.
			\end{itemize}
			In practice, the two methods are both very useful, but they have individual application scenarios.
		\end{remark}
		\jpang{
			\begin{remark}
				Lyapunov optimization and OCO have two things in common in that they both do not require historical data and prediction. 
				Other online optimization methods for multi-period problems either need historic data or prediction, such as reinforcement learning (RL) and receding-horizon control (RHC). In RL, a large amount of data is needed to train the controller offline before it can be used online at a reasonable performance in the time-varying environment. In RHC, also known as model predictive control (MPC), predictions on the near future data are required. It is important to note that RHC proceeds in a rolling-horizon manner. As such, both of them are computationally more expensive compared with Lyapunov optimization and OCO, especially for large-scale problems with long time horizons. Although they both have fruitful research achievements, we do not include them in this review due to data and prediction limitations.
		\end{remark}}
		
		\section{Challenges and Prospects}\label{Future_direction}
		
		This section presents the critical challenges and several future directions of online optimization in power systems, i.e., the capability of plug-and-play, transient performance enhancement, and online pursuit of Nash equilibrium (NE).
		
		\subsection{Capability of plug-and-play}
		
		With the proliferation of DGs, the power systems will no longer be dominated by a few large SGs, but by a massive amount of small devices with various dynamical characteristics. Because they belong to individuals instead of the utility company, some DGs may switch off or switch on unexpectedly. This requires the controller to have the capability of what we term plug-and-play, which includes two stages: 1) stability guarantees after plugging in; 2) achieving optimality in playing. Both of these should be realized in an online fashion. 
		
		\subsubsection{Online stability guarantee}
		Stability is the primary concern of power systems, which is usually analyzed by offline methods, such as time-domain simulation, eigenvalue analysis, and direct methods based on the energy function. {\color{black}The first and third approaches require a detailed system model. The second one is applicable to a linear system model with small disturbances. To summarize, all of the existing methods require the complete model of the system, which is available in the traditional power system consisting of several large SGs. However, it is difficult to realize due to the unexpected (dis)connection of DGs under the new situation. In addition, these methods rely on a known equilibrium point, which are also less applicable due to the uncertainty and volatility of renewable generation sources. }
		Thus, it is necessary to develop an online stability analysis which is adaptable to the volatile environment. Recently, incremental passivity theory has attracted much attention, as it eliminates the influence of equilibrium fluctuation. Some works derive local stability criteria, which provide verifiable conditions for DGs to connect \cite{Peng2020Distributed_stability}. If DGs satisfy this condition, they can be integrated. Otherwise, it will be not allowed to connect. Because the condition can be checked by local variables, it is easy to use and suitable for online applications. {\color{black}Although it is inspiring, much work still needs to be done to form a complete theoretical framework. First, the current condition is conservative, and often the stability criterion has to be made mild as it is preferred to enlarge the stability region. Second, because it is impractical to change the controllers of many already installed DGs, designing supplementary control strategies is very important to drive DGs towards satisfying these conditions. In addition, stability theory addressing varying dimensions should be developed. 
		}
		
		\subsubsection{Online optimality guarantee}
		When a DG is switched on or switched off, the optimization problem should be changed accordingly. The traditional centralized decision-making paradigm is subject to performance limitations in the situation of unexpected (dis)connection of massive DGs, such as a single point of failure, limited flexibility, and scalability, which is inapplicable for online implementation. Recently, the prevalence of distributed optimization alleviates this problem, which naturally renewed interest in this area \cite{yang2019survey}. Many distributed algorithms are developed, such as consensus-based methods \cite{Olfati:Consensus}, dual decomposition \cite{palomar2006tutorial}, ADMM \cite{boyd2011distributed}, etc. 
		For these distributed algorithms, a basic assumption should hold, i.e., each agent is equipped with an ideal solver or iterative formula. Moreover, such solvers or formulae are usually identical and highly problem-oriented. 
		Although such settings eliminate some drawbacks of the centralized paradigm, they still need to be improved in terms of adaptability. For instance, some agents may be reluctant to share their sensitive information with other agents or even a third party. Then, some works develop distributed algorithms with arbitrary local solvers  \cite{smith2018cocoa, zhang2021protocol}, which allow each agent to perform its computation locally through individual solvers, the so-called arbitrary solver. These solvers are self-customized, which increases the capability of plug-and-play. {\color{black}However, much work remains to be needed in this area. For example, asynchrony should be considered between arbitrary solvers to enable agents to operate at different paces. Moreover, an highly efficient communication topology should be designed for the online optimality pursuit.}

		\subsection{Online transient performance enhancement}
		{\color{black}OGDC is designed to achieve optimality in the steady state. 
			The transient performance determines how to reach the optimal steady state, which is also very important in dynamic control. However, most of the existing works pay less attention to the transient process. For example, in the optimal frequency control, the frequency nadir/overshoot and recovery time are critical for the system stability, which also need to be optimized besides the steady state.} To enhance the transient performance, model-based and model-free methods are both investigated. 
		In the first category, some inspiring works study the influencing factors of the transient index by spectral approach \cite{guo2017spectral}, which are derived only with linear models. MPC is also used in the frequency control problem with time-coupled state variables and constraints \cite{Ademola2021Frequency}, which is capable to compute optimal control commands based on predictions of future states and disturbances. Due to the strong uncertainties, accurate predictions or forecasts are often not possible. Moreover, the system is too complex to get an analytical model with the high penetration of renewable generations. These restrictions limit the wide application of model-based approaches. 
		For the model-free methods, RL has attracted surging attention \cite{yang2020reinforcement,chen2022reinforcement}, which makes decisions based on the data instead of explicit models. The data-driven nature allows it to adapt to uncertain dynamical environments by incremental learning \cite{wang2022instance}. {\color{black}Thus, RL has the potential of outperforming model-based methods in transient performance enhancement if the detailed model cannot be obtained.} 
		
		Recently, an interesting research topic on integrating model-based and model-free methods has emerged. Its motivations are twofold: 1) although the power system model is not accurate, the model-based methods still have acceptable performance in the application; 2) purely model-free approaches suffer from inherent limitations, such as scalability, sample complexity, and the heavy computation burden. Combining model-based and model-free methods together may achieve the benefits of both \cite{qu2020combining}. For instance, in OGDC, the current model-based controllers are designed for optimal operation in the steady state, while a model-free supplementary model could be added to enhance the transient performance. The model-based approach can provide a warm start for the model-free controllers. Despite limited works on this subject so far, combining model-free with model-based methods is envisioned to be a promising direction.
		{\color{black}Again, much work is left to be done in the future. First, the integration structure of model-free and model-based controllers should be well designed. Potential ways include implementing two controllers in serial, in parallel, or embedding one as an inner module in the other \cite{Qi2019Integrating}. Generally speaking, parallel integration does not break the original control structure and seems to be the easiest way in practice. In addition, the interpretability should be strengthened. Understanding the mechanism of many model-free controllers based on neural networks or machine learning remains an open question. With the help of the model-based method, its interpretability is improved tremendously. Then, the transient performance could be guaranteed theoretically.
		}
		
		{\color{black}
			\subsection{Online optimization with predictions}
			In the OCO, no future information is needed. With the advancement of techniques, short-term predictions of renewable generations and loads may be made available. The algorithm performance will be greatly improved if future information is properly utilized. In this situation, at each stage $t$, an agent could observe uncertainties for the next $W$ stages, and then makes a decision $x_{t}$. Generally, MPC naturally fits such problems, which require solving $ W $-stage optimization problems at each stage. The OCO with predictions avoids time-consuming computation, which is much faster. 
			This variant has attracted a lot of attention in recent years. Some algorithms have been proposed and their optimality guarantees have been analyzed, which shows that the performance is improved greatly, e.g., the regret lower bound decays exponentially with $W$ \cite{Yingying2021Online}. One potential direction is to consider the influence of the inaccuracy of predictions. If the prediction noises are modeled as uncertainty sets, the robustness of the algorithm should be analyzed. Similarly, if we know their distributions, the confidence interval of the regret can be given. Moreover, temporal-coupled constraints widely exist in power systems, such as the SoC update of the energy storage system and ramping limits, which are seldom considered in the OCO framework. To solve such problems, the key difficulty is to find temporal decomposition methods, such as dynamic programming. This also motivates another research direction, i.e., integrating reinforcement learning with OCO. With future predictions, the learning efficiency may be improved greatly by rolling out. 
			
		}
		
		\subsection{Online pursuit of generalized Nash equilibrium}
		
		With the large penetration of DGs and ESSs, the distribution network has been witnessing the emergence of massive ``prosumers", i.e., proactive consumers, which can both produce and consume electric power \cite{Morstyn2018Using}. In addition, the advancement of communication and control technologies at the consumer level enables them to communicate with each other. These changes encourage energy trading among prosumers, which forms local ``Peer-to-Peer" (P2P) power markets \cite{wang2021distributedGNE}. In this market, each prosumer attempts to maximize his profit while satisfying physical laws. 
		The market will be cleared until NE is reached, from which no player can unilaterally deviate for a higher profit \cite{Jianrui2022Cooperative}. Considering the global constraints, the power markets in distribution networks are formulated by a generalized Nash game (GNG), while the corresponding NE is also called generalized Nash equilibrium (GNE). 
		Due to the strong volatility and fast fluctuation, the GNE of these markets also needs to be pursued online. Similar to online optimization, the potential topics of the online game also include three aspects: NE-guided dynamic control (NEGDC), feedback GNE seeking, and online multi-stage GNE seeking.

		\subsubsection{NE-guided dynamic control}
		{\color{black}The uncertainty and volatility of renewable energies and loads make the system state change rapidly. As such, it is difficult for the market to respond to the fast change of power imbalance in real-time, which seriously affects the stability and optimality of the system. Thus, the market dynamics should be combined with the inverter control of the prosumer, which leads to the NEGDC.}
		Different from the OGDC, the NEGDC is to design a (dynamic) feedback controller for each agent, which steers the system states to the GNE of the market game. 
		Reference \cite{de2019feedback} designs controllers for generators to maximize their profit by regulating the power output according to the Cournot competition. Other related works focus on the GNE-seeking dynamics, i.e., solving the game by continuous-time algorithms \cite{zhang2020distributed,zhu2021generalized,deng2022distributed}. Although the above works do not consider complicated constraints in power systems, they are still inspiring for NEGDC in controller design and stability analysis.
		Despite limited works on this subject so far, designing dynamic controllers to steer the system to the GNE of the P2P market is envisioned to be a promising direction in the power system with massive prosumers.

		\subsubsection{NE-seeking Feedback}
		The NE-seeking feedback is an extension of feedback optimization methods to a game-theoretic setting. 
		Different from recent advances in the offline GNE seeking algorithms \cite{yi2017distributed}, the online measurement feedback from the physical system will be incorporated \cite{agarwalgame}. 
		Current works mainly focus on affine constraints, which still need to be further improved in many aspects to satisfy the complicated nonlinear constraints in power systems like power flow equations. By the measurement feedback, the numerical computation of these complex constraints can be replaced by physical laws. {\color{black}This not only reduces the computation burden but also avoids private information exchange to a certain extent. }
		Key challenges include controller design, existence and uniqueness analysis of the NE, as well as the convergence/stability guarantees. 
		
		\subsubsection{Online multi-stage NE seeking}
		{\color{black}Currently, the formulation of online multi-stage NE seeking is similar to that of OCO. Each player aims to minimize its own cost function subject to some global constraints, i.e. $\min\limits_{x_{i, t} \in \Omega_{i}} \ J_{i, t}\left(x_{i, t}, x_{-i, t}\right), \ s.t.\  g_{t}\left(x_{t}\right) \leq \mathbf{0} $, where both the objectives and constraints could be time-varying \cite{meng2021decentralized,Kaihong2021Online}. Similar to the OCO, regret can also be used to evaluate the performance of the algorithm. Current problem formulation does not consider the coupling of adjacent time-steps, which, however, widely exists in power systems, such as the operation of ESSs. If we consider that, the problem adopts the following form $\min \limits_{x_{i, t} \in \Omega_{i}} \ \sum\limits_{t=1}^T J_{i, t}\left(x_{i, t}, x_{-i, t}\right), \ s.t.\  g_{t}\left(x_{t}, x_{t+1}\right) \leq \mathbf{0} $. The temporal-coupled constraint is not a trivial extension, which will bring more challenges to the analysis framework, including the NE existence and algorithm design. Lyapunov optimization approach could be used to decouple the time coupling, but this needs further investigation on whether it can help solve a GNG. Another extension then would be how to modify it to multi-player setting and how to evaluate the deviation of the result from the GNE. }
		
		{\color{black}To summarize, the online pursuit of GNE is more challenging compared with online optimization. First, the existence and uniqueness of GNE are difficult to justify. The variational inequality approach is utilized to find a GNE, which can guarantee the existence of GNE with assumptions on the monotonicity of the pseudo-gradient. However, this still cannot guarantee the uniqueness \cite{yi2017distributed}. Second, information sharing is constrained. Because each prosumer is an independent stakeholder to maximize their own profit in a competitive market, they may be reluctant to share private information. Thus, the privacy preservation mechanism should be well designed. Third, the performance analysis is more challenging, since they include multiple metrics such as stability and regret. The pseudo-gradient is widely adopted in the algorithm design, which, however, has no symmetry compared with the gradient in online optimization. This implies that the pseudo-Hessian matrix is not symmetric. Consequently, many second-order properties do not exist, and cannot be applied in the performance analysis.}

		
		\section{Conclusion}\label{Conclusion}
		
		Although some works have been devoted to the online optimization of power systems with high penetration of renewable generations, they have different interpretations based on their corresponding time scales, which leads to much confusion. In this paper, we provide a comprehensive review and comparative analysis of three different online optimization notions in power systems, including motivations, time scales, popular algorithms, theoretic foundations, and typical applications.  
		{\color{black}Moreover, we also present the critical challenges and several future directions, such as the capability of plug-and-play, transient performance enhancement, online optimization with predictions, and online pursuit of NE. }
		\jpang{It must be emphasized that online optimization is strongly problem-dependent, and has no unified mathematic paradigm.
			Our hope is that this paper helps the readers to figure out where and when to use what types of online optimization algorithms.}

		\bibliographystyle{IEEEtran}
		\bibliography{mybib}
		

		\ifCLASSOPTIONcaptionsoff
		\newpage
		\fi

	\end{document}